
\documentstyle{amsppt}
\magnification 1200
\def\boxit#1{\vbox{\hrule\hbox{\vrule\kern2.5pt\vbox{\kern2.5pt#1
\kern2.5pt}\kern2.5pt\vrule}\hrule}}
\nologo
\NoBlackBoxes
\input epsf.tex
\topmatter
\title Normal Presentation on Elliptic Ruled Surfaces \endtitle
\author Francisco Javier Gallego \\ and \\ B. P. Purnaprajna
\endauthor
\affil Brandeis University \\ Department of Mathematics \endaffil
\address{Francisco Javier Gallego: Departamento de Algebra,
Facultad de Matematicas, Uni\-ver\-si\-dad Com\-plu\-ten\-se de
Madrid, 28040 Madrid, Spain }\endaddress
\email{ gallego\@sunal1.mat.ucm.es}\endemail
\address{ Bangere P. Purnaprajna:
Dept of Mathematics
 Brandeis University,
Waltham MA 02254-9110}\endaddress
\email{ purna\@max.math.brandeis.edu}\endemail
\thanks{We would like to thank our advisor David Eisenbud for
his encouragement and helpful advice. We are also glad to
thank Aaron Bertram, Raquel Mallavibarrena and Giuseppe Pareschi
for helpful
discussions.}\endthanks
\endtopmatter
\document
\headline={\ifodd\pageno\rightheadline \else\leftheadline\fi}
\def\rightheadline{\tenrm\hfil
\eightpoint NORMAL PRESENTATION ON ELLIPTIC RULED SURFACES
 \hfil\folio}
\def\leftheadline{\tenrm\folio\hfil
\eightpoint F.J. GALLEGO \& B.P. PURNAPRAJNA \hfil}
\vskip .3 cm

\roster
\item ""{\bf  Introduction}
\item "{\bf 1.}" {\bf Background material}
\item "{\bf 2.}" {\bf General results on normal presentation}
\item "{\bf 3.}" {\bf Ampleness, base-point-freeness and cohomology of
line bundles
on elliptic ruled surfaces}
\item "{\bf 4.}" {\bf Normal presentation on elliptic ruled
surfaces}
\item "{\bf 5.}" {\bf Koszul algebras}
\item "" {\bf References}
\endroster
\heading Introduction \endheading

This article deals with the {\it normal presentation} of line bundles over an
elliptic ruled surface.
 Let $X$ be an irreducible projective
variety and $L$ a very ample line bundle on $X$,  whose complete linear
series defines $$\phi_L : X \longrightarrow \bold P (\text H^0(L))$$
 Let $S=\bigoplus_{n=0}^\infty \text S^m\text H^0(X,L) $ and  Let  $ R (L)
=
\bigoplus_{n=0}^\infty\text H^0(X, L^{\otimes n})$  be the homogenous
coordinate ring associated to $L$.  Then
$ R$ is a finitely generated  graded module over S, so it has a minimal graded
 free resolution. We say  that  the line bundle $L$ is normally generated
if the natural maps
$$ \text S^m\text H^0(X,L) \to \text H^0(X,L^{\otimes m}) $$
are surjective for all $m \geq 2$. If $L$ is normally generated, then we say
 that  $L$ satisfies property $N_p$,  if the matrices in the free
resolution of
$R$ over $S$ have linear entries until the  pth  stage. In particular,
property
$N_1$ says that the homogeneous ideal $I$ of $X$ in $\bold P (\text H^0(L))$ is
generated by quadrics. A line bundle satisfying property $N_1$ is also called
normally presented.

 Let $R= \bold k \oplus R_1 \oplus R_2
\oplus
\dots$ be a graded  algebra over a field  $\bold k$. The algebra
$R$ is a Koszul ring iff Tor$_i^R(\bold k,\bold k)$
has pure degree $i$ for all $i$.

In this article we determine exactly (Theorem 4.1) which line bundles on
elliptic ruled surface $X$
 are normally presented (Yuko Homma has classified in [Ho1] and [Ho2] all line
bundles which are normally generated on an elliptic ruled surface). In
particular we see that numerical classes of normally presented divisors  form a
convex set. (See Figure  1 for the case $e(X)=-1$; recall that $\text{Num}(X)$
is generated by the class
 of a minimal section $C_0$ and by the class of a fiber $f$ and that $C_0$ is
ample.) As a corollary of the above result we show that {\it Mukai's
conjecture} is true
 for the normal presentation of the {\it adjoint linear series} for an elliptic
ruled surface.

In section 5 of this article, we show that if $L$ is normally presented on $X$
then the homogeneous
coordinate ring associated to $L$ is Koszul. We also give a new proof of the
following
result due to Butler:
 if deg$(L) \geq 2g+2$ on a curve $X$ of genus g, then $L$ embeds $X$ with
Koszul
homogeneous coordinate ring.

 To put things in perspective, we would like to recall what is known regarding
these
questions in the case of curves. A classical result of
Castelnuovo  (c.f. [C]) says
 that if deg$(L) \geq 2g+1$,
 $L$ is
normally  generated. St. Donat and Fujita ([F] and [S-D]) proved that if
deg$(L)
\geq 2g+2$, then $L$ is normally presented.
These theorems have been recently
generalized to higher syzygies by Green (see [G]),
who proved that if deg$(L) \geq 2g + p +1$, then $L$ satisfies the
property $N_p$. One way of generalizing the above results to higher dimensions
is to interpret them in terms of adjoint linear series: let
$\omega _X$  be the canonical bundle of a curve $X$, and let $A$ be an ample
line bundle (since
$X$ is a curve,
$A$ is ample iff deg$(A) > 0$). If
$L = \omega _X
\otimes A ^{\otimes 3}$ (respectively $L = \omega _X
\otimes A ^{\otimes p+3}$), then Castelnuovo's Theorem
(respectively Green's Theorem) says that $L$ is normally
generated  (respectively satisfies property $N_p$).

Unlike the case of curves, the
landscape of surfaces (not to speak of higher dimensions) is
relatively uncharted.  Recently  Reider proved (c.f.
[R]) that   if
$X$ is a surface over the complex numbers, then
$\omega _X
\otimes A ^{\otimes 4}$ is very ample. Mukai has conjectured that $\omega _X
\otimes A ^{\otimes p + 4}$ satisfies $N_p$. Some work in this direction
has been done by David Butler in [B], where he studies the syzygies
of adjoint linear series  on ruled varieties. He proves
that if the dimension of $X$ is $n$, then $\omega _X
\otimes A ^{\otimes 2n+1}$ is normally generated and $\omega _X
\otimes A ^{\otimes 2n+2np}$ is normally presented;
specializing to the case of ruled surfaces, his result says
that $\omega _X
\otimes A ^{\otimes 5}$ is normally generated and that $\omega _X
\otimes A ^{\otimes 8}$ is normally presented. In this article we consider not
just
the adjunction bundle, but any very ample line bundle
on an elliptic ruled surface.
 In particular, we prove that $\omega _X
\otimes A ^{\otimes 5}$ is normally presented thereby proving Mukai's
conjecture for $p=1$ in the case of elliptic ruled surface.

In a sequel to this article we generalize our results on normal presentation
to higher syzygies.  We there show the following:  let $L=B_1\otimes
...\otimes B_{p+1}$ be a line bundle on
$X$, where each $ B_i$ is base point free and ample, then $L$ satisfies
property
$N_P$. As a corollary we show that $\omega _X
\otimes A ^{\otimes 2p+3}$ satisfies property $N_p$.

\heading  1. Background material \endheading
{\bf Convention.} Throughout this paper we work
over an algebraic closed field $\bold k$.
\vskip .2 cm
We state in this section some results we will use later. The first
one is this beautiful cohomological characterization by Green of the
property $N_p$. Let $L$ be a globally generated line bundle. We
define   the vector bundle
$M_L$ as follows:
$$0 \to M_L \to \text H^0(L) \otimes \Cal O_X \to L \to 0 \ .\eqno
(1.1)$$
In fact, the exact sequence (1.1) makes sense for any
variety
$X$ and any vector bundle $L$ as long as $L$ is globally generated.

\proclaim {Lemma 1.2} Let $L$ be a normally generated line bundle
on a variety $X$ such that $\text H^i(L^{\otimes 2-i})=0$ for all
$i \geq 1$. Then,
$L$ satisfies the property
$N_p$ iff
\text{$\text H ^1(\bigwedge ^{p'+1} M_L \otimes L)$} vanishes for all $1 \leq
p' \leq p$.
\endproclaim
{\it Proof.} The lemma is a corollary of [GL], Lemma 1.10.
\boxit{}

\vskip .2 cm

 (1.2.1) If the $\text{char}(\bold k)\neq 2$, we can obtain the
vanishing of  $\text H ^1(\bigwedge ^{2} M_L\otimes L)$  by showing
the vanishings of
$\text H ^1(M_L^{2}\otimes L)$, because
$\bigwedge ^{2} M_L\otimes L$ is in this case a direct summand of
$M_L^{2}\otimes L$.

\vskip .2 cm

The other main tool we will use is a generalization by Mumford
of a lemma of Castelnuovo:
\proclaim {Theorem  1.3} Let $L$ be a base-point-free line bundle
on a variety $X$ and let $\Cal F$ be a coherent sheaf on $X$. If
$\text H^i (\Cal F \otimes L ^{-i}) = 0$ for all $i \geq 1$, then
the multiplication map
$$ \text H ^0(\Cal F \otimes L ^{\otimes i}) \ \otimes \text
H ^0(L) \to  \text H ^0(\Cal F \otimes L ^{\otimes i+1})$$
is surjective for all $i \geq 0$.
\endproclaim

{\it Proof.} [Mu], p. 41, Theorem 2. Note that the assumption made
there of $L$ being ample is unnecessary. \boxit{}

\vskip .2 cm
It will be useful to have the following characterization of
projective normality:
\proclaim {Lemma 1.4} Let $X$ be a surface with geometric genus
$\text h^2(\Cal O_X)=0$ and let
$L$ be an ample, base-point-free line bundle.
If $\text H^1(L )
=0$, then $L$ is normally generated iff $\text H^1(M_L \otimes
L) =0$.
\endproclaim

{\it Proof.}
The line bundle $L$ is normally generated iff the map
$$ \text S^m\text H^0(X,L) @>\alpha>> \text H^0(X,L^{\otimes m})$$
is surjective.
The map $\alpha$ fits in the following commutative diagram:
$$ \matrix
\text H^0(L)^{\otimes m} &@>\beta>> &\text S^m\text H^0(L)\cr
@VV\gamma_1V&\cr
\text
H^0(L^{\otimes 2}) \otimes
\text H^0(L)^{\otimes m-2}&\cr
@VV\gamma_2 V @VV \alpha V \cr
\vdots&\cr
@VV\gamma_{m-2}V&\cr
\text H^0(L^{\otimes m-1})
\otimes H^0(L)&@>\gamma_{m-1}>>&\text H^0(L^{\otimes m}) \cr
\endmatrix
$$

The map $\beta$ is surjective. From this fact it follows that the
surjectivity of
$\alpha$ is equivalent to the surjectivity of
\text{$\gamma_{m-1}
\circ \cdots \circ
\gamma_1$.} Theorem 1.3 implies the
surjectivity of
$\gamma_2,
\dots,
\gamma_{m-1}$. Hence the surjectivity of $\gamma_1$ implies the
surjectivity of $\alpha$. On the other hand, if $m=2$ the
surjectivity of $\alpha$ implies the surjectivity of $\gamma_1$.
Finally from (1.1) we obtain
$$\text H^0(L)\otimes\text H^0(L) \to \text H^0(L^{\otimes 2}) \to
\text H^1(M_L \otimes L) \to \text H^0(L) \otimes \text H^1(L) \ .$$
Therefore the vanishing of $\text H^1(L)$ implies that the
surjectivity of
$\gamma_1$ is equivalent to the vanishing of $\text H^1(M_L \otimes L)$.
\boxit{}

\heading 2. General results on normal presentation
\endheading

As mentioned in the introduction, according to  our
philosophy, the tensor product of two
base-point-free line bundles $B_1$ and $B_2$ (provided it
is ample and that certain higher cohomology groups vanish) should
be normally presented. This philosophy is made concrete in the
following

\proclaim {Proposition 2.1}
Let $X $ be a surface with geometric genus $0$ and $B_1$ and
$B_2$ base-point-free line bundles such that
\text{$\text H^1(B_1)=\text H^1(B_2)=\text H^2(B_2
\otimes B_1 ^*)=$}$\text H^2(B_1 \otimes B_2 ^*)=0$ and let $L= B_1
\otimes B_2$. Then $\text H ^1(M_L^{\otimes p+1}\otimes
L)=0$ for $p = 0,1$. In particular, if
$L$ is ample and ${\text{\rm{char}}}(\bold k) \neq 2$, then $L$ is
normally presented.
\endproclaim

We will prove a more general version of Proposition 2.1 in
 section 5.
\vskip .2 cm
{}From this proposition
we will obtain corollaries for Enriques surfaces (Corollary 2.8) and
for elliptic ruled surfaces (Theorem  4.1). To prove Proposition 2.1
we will need several lemmas and observations:
\proclaim
{Observation 2.2} Let $X$ be a surface with geometric genus $0$,
let $P$ be an
effective line bundle, and let $B$ be a  line bundle such
that, for some $\frak p \in
|P|$, $B \otimes \Cal O_{\frak p}$ is trivial or has a global section
vanishing at finite subscheme of $\frak p$  (e.g., let $B$ be
base-point-free). If $\text H^1(P)=
\text H^1(B) =0$, then $\text H^1(B\otimes P)=0$.
\endproclaim


\proclaim
{Observation 2.3} Let $X$ be a surface, let $P$ be an
effective line bundle and L a coherent sheaf. If H$^2(L) =0$, then
H$^2(L \otimes P) =0$
\endproclaim

\proclaim {Lemma 2.4}
Let $X$ be a surface, let $B$ be a globally generated line bundle such
that H$^1(B) = 0$ and
let $Y$ be a curve in $X$ such that $B \otimes \Cal O_Y (Y) $ is globally
generated. Then $B \otimes \Cal O_X (Y)$ is also globally
generated.
\endproclaim

{\it Proof.}
The result is the surjectivity of the middle vertical arrow
in the following commutative diagram:
$$\matrix
\text H ^0(B) \otimes \Cal O _X &\hookrightarrow &\text H ^0(B
\otimes \Cal O _X (Y)) \otimes \Cal O _X&\twoheadrightarrow
&\text H ^0(B
\otimes \Cal O _Y (Y)) \otimes \Cal O _X\\
\downarrow && \downarrow
&&\downarrow
\\ B&\hookrightarrow&B
\otimes \Cal O _X (Y)&\twoheadrightarrow &B
\otimes \Cal O _Y (Y)\ .\\
\endmatrix$$
The hypothesis is that the vertical left hand side arrow and the
vertical right hand side arrow are surjective. \boxit{}
\vskip .5 cm

\proclaim {Lemma 2.5}
Let $X $ be a surface with geometric genus $0$, let $B_1$ and
$B_2$  be two base-point-free line bundles and let $L= B_1
\otimes B_2$.  If
$\text H^1(B_1)=\text H^1(B_2)=0$
 and $\text H^2(B_2 \otimes B_1 ^*) =0$,
    then
$\text H^1(M_L \otimes B_1 ^{\otimes n}) = 0$ for all $n \geq 1 $
\endproclaim
{\it Proof.}
If we tensor  exact sequence (1.1)
with
$B_1^{\otimes n}$ and take global sections, we obtain
$$\displaylines{ \text H ^0(L) \otimes \text H ^0(B_1 ^{\otimes n})
@>
\alpha >>
\text H ^0(L
\otimes B_1^{\otimes n})\cr
 \to \text H ^1(M_L \otimes B_1^{\otimes
n}) \to \text H ^0(L)
\otimes \text H ^1(B_1^{\otimes n})\ .}$$
{}From Observation 2.2
it follows that the vanishing of H$^1(M_L \otimes B_1^{\otimes n})
$ is equivalent to the surjectivity of the multiplication map
$\alpha$. In the case $n=1$ the surjectivity of $\alpha$ follows
trivially from our hypothesis and Theorem 1.3.
The proof of the  surjectivity
of
$\alpha$ goes by induction. We show here only the case $n=2$.  We
consider the  commutative diagram
$$\matrix
\text H ^0(B_1) \otimes \text H ^0(B_1) \otimes \text H
^0(L)&\to&\text H ^0(B_1^{\otimes 2}) \otimes \text H ^0(L)\\
@VV\gamma V @VV\alpha V\\
\text H^0(B_1) \otimes \text H
^0( B_1 \otimes L)&@> \delta >>&\text
H ^0(B_1 ^{\otimes 2}
\otimes L)\\
\endmatrix$$
where the maps are the obvious ones coming from multiplication. To
prove the surjectivity of $\alpha $ it suffices to prove that
$\gamma$ and
$\delta$ are surjective. The surjectivity of $\gamma$ follows
from the surjectivity of $\alpha$ when $n = 1$. To prove
the surjectivity of
$\delta$, again  by Theorem 1.3, it is enough to check that
$\text H ^1(L) =
\text H ^2(B_2 ) = 0$. This follows from the hypothesis and from
the Observations 2.2 and 2.3. \boxit{}

\proclaim {Lemma 2.6}
Let $X $ be a surface with geometric genus $0$, let $B_1$ and
$B_2$  be two base-point-free line bundles and let $L= B_1
\otimes B_2$ be nonspecial. Let $B_1$ and $B_2$ satisfy the
conditions   $\text H^1(B_1^{\otimes 2})=\text H^1(B_2)=0$
 and $\text H^2(B_2 \otimes B_1 ^*) =\text H^2(B_1^{\otimes 2}
 \otimes B_2 ^*)=0$. If $P$ is any effective  line bundle
 on
$X$ such that either H$^1(P)=0$ or $P \simeq \Cal
O$,   then
H$^1(M_L \otimes L \otimes P) = 0$.
In particular, if $L$ is ample, $L$ is normally generated.
\endproclaim

{\it Proof. }
If we tensor (1.1) with $L\otimes P$ and take global
sections, we obtain
$$ \text H ^0(L) \otimes \text H ^0(L \otimes P) @> \alpha >> \text
H ^0(L ^{\otimes 2}
\otimes P) \to \text H ^1(M_L \otimes L \otimes P) \to \text H
^0(L)
\otimes \text H ^1(L \otimes P) \ .$$
{}From Observation 2.2 it follows that  $\text H ^1(L \otimes
P)$ vanishes. Therefore the vanishing of
\text{$\text H ^1(M_L
\otimes L \otimes P)$} is equivalent to the surjectivity of the
multiplication map $\alpha$. To prove the surjectivity of
$\alpha$ we use the same trick as in the proof of the previous lemma. We
write this
 commutative diagram:
$$\matrix
\text H ^0(B_2) \otimes \text H ^0(B_1) \otimes \text H
^0(L\otimes P)&\to&\text H ^0(L) \otimes \text H ^0(L \otimes P)\\
@VV\gamma V @VV\alpha V\\
\text H^0(B_2) \otimes \text H
^0(L\otimes B_1 \otimes P)&@> \delta >>&\text
H ^0(L ^{\otimes 2}
\otimes P) \ .\\
\endmatrix$$
 It
suffices then to prove that $\gamma$ and $\delta$ are surjective
and by Theorem 1.3 it is enough to check that \text{$\text H ^1(B_2
\otimes P) =
\text H ^2(B_2 \otimes B_1^* \otimes P) = \text H
^1(B_1^{\otimes 2} \otimes P) =$}$ \text H ^2(B_1^{\otimes 2}
\otimes B_2^*\otimes P)=0$. These vanishings follow trivially
from the hypothesis of the lemma and from Observations 2.2 and
2.3. \boxit{}

\vskip .4 cm
(2.7)
{\it Proof of Proposition 2.1.}  Observation 2.2  implies that
$\text H ^1(L)$ vanishes. Thus from Lemma 2.6 it follows that $\text
H ^1(M_L
\otimes L )=0$. This implies that the vanishing of  \text{$\text H
^1(M_L ^{\otimes 2}
\otimes L )$} is equivalent to the surjectivity of
the multiplication map
$$\text H ^0(M_L \otimes L
 )
\otimes H ^0(L) @> \alpha >> \text H ^0(M_L \otimes L ^{\otimes
2} )\ . \eqno(2.7.1)$$ To prove the surjectivity of $\alpha$ we write
this commutative diagram:
$$\matrix
\text H ^0(B_2) \otimes \text H ^0(B_1) \otimes \text H
^0(M_L\otimes L)&\to&\text H ^0(L) \otimes \text H ^0(M_L \otimes
L)\\ @VV V @VV \alpha V\\
\text H^0(B_2) \otimes \text H
^0(M_L \otimes L\otimes B_1 )&\to&\text
H ^0(M_L \otimes L ^{\otimes 2}
)\ .\\
\endmatrix$$
 By Theorem 1.3 it is enough to check that
$$ \displaylines{\text H ^1(M_L \otimes B_2) =
\text H ^1(M_L \otimes B_1 ^{\otimes 2}) =0 \cr
\text H ^2(M_L \otimes B_2 \otimes B_1 ^*)=
\text H^2(M_L \otimes B_1 ^{\otimes 2} \otimes B_2 ^*)=0 \ .\cr}
$$

The first two vanishings follow from Lemma 2.5. The other two follow from
sequence (1.1) and from Observations 2.2
and 2.3.

If $L$ is ample, it follows from Lemma 1.2, (1.2.1) and Lemma 1.4 that $L$ is
normally presented. \boxit{}

\vskip .2 cm
The conditions on the vanishing of cohomology required in the
statement of Proposition 2.1 are not so restrictive. For instance if
we take $B_1$ and $B_2$ equal and ample, the conditions on the
vanishing of
$\text H^2$ are automatically satisfied for surfaces with
geometric genus $0$. If the surface we are considering is Enriques
or elliptic ruled the vanishing of $\text H^1$ also occurs. The
next corollary is an outcome of these observations.

\proclaim {Corollary 2.8 }
Let $X$ be an Enriques surface, let $\text{\rm{char}}(\bold k) =0$
and let
$B$ be an ample line bundle on $X$ without base points. Then
$B^{\otimes 2}$ is normally presented.
\endproclaim

{\it Proof.}
Since $K_X \equiv 0$ and $B$ is ample, $\omega_X \otimes B$ is also
ample and by Kodaira vanishing, H$^1(B)=0$. Thus we can apply
Proposition 2.1. \boxit{}

\vskip .2 cm

\heading 3. Ampleness, base-point-freeness and cohomology of line bundles
 on elliptic ruled surfaces \endheading

We have shown in Corollary 2.8 that $B^{\otimes 2}$ is
normally presented if $B$ is an ample, base-pont-free line bundle
over an Enriques surface. The same result is true in the case of
elliptic ruled surfaces. However in this case we can do much
better. In fact we will be able to characterize (c.f. Theorem 4.2) those
line bundles
which are normally presented. From the statement of
Lemma 2.6 it is clear that the knowledge of the vanishing of higher
cohomology of line bundles on elliptic ruled surfaces will be
crucial for this purpose. On the other hand once we know that the
tensor product of two base-point-free line bundles is normally
presented, knowing in addition which line bundles on an elliptic
ruled surface are base-point-free will allow us to characterize
those line bundles that are normally presented. In this light we
will devote this section  to recalling the vanishing of
cohomology of line bundles and  the characterization of
base-point-free line bundles on elliptic ruled surfaces, .

We introduce now some notation and recall some
elementary facts about elliptic ruled  surfaces. Proofs
for the statements of this paragraph can be found in [H], \S V.2. In
this and the next section
$X$ will denote a smooth elliptic ruled surface, i.e. $X = \bold P
(\Cal E)$, where $\Cal E$ is a vector bundle of rank 2 over a
smooth elliptic curve $C$.  We will assume $\Cal E$ to be
normalized, i.e., $\Cal E$ has global sections but twists of it
by  line bundles of negative degree do not. Let
$\pi$ denote the projection from
$X$ to
$C$.  We set $\Cal O(\frak e )=
\bigwedge ^2 \Cal E$ and
$e = -\text{deg} \ \frak e \geq -1$. We fix a minimal section $C_0$
such that
$\Cal O (C_0) = \Cal O_{\bold P (\Cal E)}(1)$. The group Num$(X)$ is
generated by
$C_0$ and by the class of a fiber, which we will denote by $f$. If
$\frak a$ is a divisor on $C$, $\frak a f$ will denote the
pullback of $\frak a$ to $X$ by the projection from $X$ to $C$.
Sometimes,
when $\text{deg}\ \frak a =1$, we will write, by an abuse of
notation,
$f$  instead of
$\frak a f$. The
canonical divisor
$K_X$ is linearly equivalent to $-2C_0 + \frak e f$, and hence
numerically equivalent to $-2C_0 -ef$.
\proclaim {Proposition 3.1}

Let $L$ be a line bundle on $X$, numerically equivalent to
$aC_0+bf$.

\quad If $e =-1$:
\vskip .1 cm
$$\vbox {\offinterlineskip \hrule
\halign{\vrule
#&\hfil\quad#\quad\hfil&\vrule#&\hfil\quad#\quad\hfil&\vrule#&
\hfil\quad#\quad\hfil&\vrule#&\hfil\quad#\quad\hfil&\vrule#&
\hfil\quad#\quad\hfil&\vrule#
\cr height .1cm &&height .1cm &&height .1cm &&height .1cm&&height
.1cm&&height
.1cm
  \cr
&$a$&&$b$&&$\text h^0(L)$&&$\text h^1(L)$&&$\text h^2(L)$&\cr
height .1cm &&height .1cm &&height .1cm &&height .1cm&&height .1cm
  &&height .1cm\cr
  \noalign{\hrule}
height .1cm &&height .1cm &&height .1cm &&height .1cm&&height
.1cm&&height .1cm
  \cr
  &&&$b > -a/2$&&$>0$&&$0$&&$0$& \cr
height .1cm &&height .1cm &&height .1cm &&height .1cm&&height .1cm
&&height .1cm  \cr
&&&\omit\hrulefill&&\omit\hrulefill&&\omit\hrulefill&
&\omit\hrulefill&\cr
height .1cm &&height .1cm &&height .1cm &&height .1cm&&height .1cm
  &&height .1cm\cr
  &$a \geq 0$&&$b = -a/2$&&$?$&&$?$&&$0$& \cr
height .1cm &&height .1cm &&height .1cm &&height .1cm&&height .1cm
  &&height .1cm \cr
&&&\omit\hrulefill&&\omit\hrulefill&&\omit\hrulefill&
&\omit\hrulefill&\cr
height .1cm &&height .1cm &&height .1cm &&height .1cm&&height .1cm
  &&height .1cm \cr
  &&&$b < -a/2$&&$0$&&$>0$&&$0$& \cr
height .1cm &&height .1cm &&height .1cm &&height .1cm&&height .1cm
  &&height .1cm\cr
\noalign{\hrule}
height .1cm &&height .1cm &&height .1cm &&height .1cm&&height .1cm
  &&height .1cm\cr
  &$a=-1$&&any $b$&&$0$&&$0$&&$0$& \cr
height .1cm &&height .1cm &&height .1cm &&height .1cm&&height .1cm
 &&height .1cm \cr
   \noalign{\hrule}
height .1cm &&height .1cm &&height .1cm &&height .1cm&&height .1cm
  &&height .1cm\cr
  &&&$b > -a/2$&&$0$&&$>0$&&$0$& \cr
height .1cm &&height .1cm &&height .1cm &&height .1cm&&height .1cm
  &&height .1cm\cr
&&&\omit\hrulefill&&\omit\hrulefill&&\omit\hrulefill&
&\omit\hrulefill&\cr
height .1cm &&height .1cm &&height .1cm &&height .1cm&&height .1cm
  &&height .1cm\cr
  &$a \leq -2$&&$b = -a/2$&&$0$&&$?$&&$?$& \cr
height .1cm &&height .1cm &&height .1cm &&height .1cm&&height .1cm
  &&height .1cm\cr
&&&\omit\hrulefill&&\omit\hrulefill&&\omit\hrulefill&
&\omit\hrulefill&\cr
height .1cm &&height .1cm &&height .1cm &&height .1cm&&height .1cm
  &&height .1cm\cr
  &&&$b < -a/2$&&$0$&&$0$&&$>0$& \cr
height .1cm &&height .1cm &&height .1cm &&height .1cm&&height .1cm
  &&height .1cm\cr}  \hrule }$$

\quad If $e \geq 0$:
\vskip .1 cm
$$\vbox {\offinterlineskip \hrule
\halign{\vrule
#&\hfil\quad#\quad\hfil&\vrule#&\hfil\quad#\quad\hfil&\vrule#&
\hfil\quad#\quad\hfil&\vrule#&\hfil\quad#\quad\hfil&\vrule# \cr
height .1cm &&height .1cm &&height .1cm &&height .1cm&&height .1cm
  \cr
&$a$&&$b$&&$\text h^0(L)$&&$\text h^2(L)$&\cr
height .1cm &&height .1cm &&height .1cm &&height .1cm&&height .1cm
  \cr
  \noalign{\hrule}
height .1cm &&height .1cm &&height .1cm &&height .1cm&&height .1cm
  \cr
  &&&$b > 0$&&$>0$&&$0$& \cr
height .1cm &&height .1cm &&height .1cm &&height .1cm&&height .1cm
  \cr    &&&\omit\hrulefill&&\omit\hrulefill&&\omit\hrulefill&\cr
height .1cm &&height .1cm &&height .1cm &&height .1cm&&height .1cm
  \cr
  &$a \geq 0$&&$b = 0$&&$?$&&$0$& \cr
height .1cm &&height .1cm &&height .1cm &&height .1cm&&height .1cm
  \cr &&&\omit\hrulefill&&\omit\hrulefill&&\omit\hrulefill&\cr
height .1cm &&height .1cm &&height .1cm &&height .1cm&&height .1cm
  \cr
  &&&$b < 0$&&$0$&&$0$& \cr
height .1cm &&height .1cm &&height .1cm &&height .1cm&&height .1cm
  \cr
\noalign{\hrule}
height .1cm &&height .1cm &&height .1cm &&height .1cm&&height .1cm
  \cr
  &$a=-1$&&any $b$&&$0$&&$0$& \cr
height .1cm &&height .1cm &&height .1cm &&height .1cm&&height .1cm
  \cr
   \noalign{\hrule}
height .1cm &&height .1cm &&height .1cm &&height .1cm&&height .1cm
  \cr
  &&&$b > -e$&&$0$&&$0$& \cr
height .1cm &&height .1cm &&height .1cm &&height .1cm&&height .1cm
  \cr    &&&\omit\hrulefill&&\omit\hrulefill&&\omit\hrulefill&\cr
height .1cm &&height .1cm &&height .1cm &&height .1cm&&height .1cm
  \cr
  &$a \leq -2$&&$b = -e$&&$0$&&$?$& \cr
height .1cm &&height .1cm &&height .1cm &&height .1cm&&height .1cm
  \cr &&&\omit\hrulefill&&\omit\hrulefill&&\omit\hrulefill&\cr
height .1cm &&height .1cm &&height .1cm &&height .1cm&&height .1cm
  \cr
  &&&$b < -e$&&$0$&&$>0$& \cr
height .1cm &&height .1cm &&height .1cm &&height .1cm&&height .1cm
  \cr}  \hrule }$$

$$\vbox {\offinterlineskip \hrule
\halign{\vrule
#&\hfil\quad#\quad\hfil&\vrule#&\hfil\quad#\quad\hfil&\vrule#&
\hfil\quad#\quad\hfil&\vrule# \cr
height .1cm &&height .1cm &&height .1cm &&height .1cm
  \cr
&$a$&&$b$&&$\text h^1(L)$&\cr
height .1cm &&height .1cm &&height .1cm &&height .1cm
  \cr
  \noalign{\hrule}
height .1cm &&height .1cm &&height .1cm &&height .1cm
  \cr
  &&&$b > ae$&&$0$& \cr
height .1cm &&height .1cm &&height .1cm &&height .1cm
  \cr    &&&\omit\hrulefill&&\omit\hrulefill&\cr
height .1cm &&height .1cm &&height .1cm &&height .1cm
  \cr
  &$a \geq 0$&&$b = ae$&&$?$& \cr
height .1cm &&height .1cm &&height .1cm &&height .1cm
  \cr &&&\omit\hrulefill&&\omit\hrulefill&\cr
height .1cm &&height .1cm &&height .1cm &&height .1cm
  \cr
  &&&$b < ae$&&$>0$& \cr
height .1cm &&height .1cm &&height .1cm &&height .1cm
  \cr
\noalign{\hrule}
height .1cm &&height .1cm &&height .1cm &&height .1cm
  \cr
  &$a=-1$&&any $b$&&$0$& \cr
height .1cm &&height .1cm &&height .1cm &&height .1cm
  \cr
   \noalign{\hrule}
height .1cm &&height .1cm &&height .1cm &&height .1cm
  \cr
  &&&$b > e(a+1)$&&$0$& \cr
height .1cm &&height .1cm &&height .1cm &&height .1cm
  \cr    &&&\omit\hrulefill&&\omit\hrulefill&\cr
height .1cm &&height .1cm &&height .1cm &&height .1cm
  \cr
  &$a \leq -2$&&$b =e(a+1)$&&$?$&
 \cr
height .1cm &&height .1cm &&height .1cm &&height .1cm
  \cr &&&\omit\hrulefill&&\omit\hrulefill&\cr
height .1cm &&height .1cm &&height .1cm &&height .1cm
  \cr
  &&&$b < e(a+1)$&&$>0$&
  \cr
height .1cm &&height .1cm &&height .1cm &&height .1cm
  \cr}  \hrule }$$

\endproclaim
\vskip .3 cm
{\it Proof.}
If $a<0$ it is obvious that h$^0(L)=0$. If $a\geq 0$, one
obtains the statements for h$^0(L)$ and h$^1(L)$ by pushing down
$L$ to $C$ and computing the cohomology there.  In
the case of $e\leq 0$, we use the fact that the symmetric powers
of $E$ are semistable bundles ([Mi], Corollary 3.7 and \S 5).
 Then we use the fact that, if $F$ is a semistable bundle over
an elliptic curve and deg$(F)>0$, then h$^0(F) >0$ and h$^1(F)
=0$ and the fact that if deg$(F)<0$, then h$^0(F) =0$ and h$^1(F)
>0$. In the case $e>0$ the computation of
cohomology on $C$ is elementary, since $E$ is decomposable. If $a=-1$,
$\pi_*L=R^1\pi_*L=0$, hence $\text H^1(L)=0$.

The other statements in the proposition follow by duality.
\boxit{}

\vskip .7 cm
The last proposition means that the vanishing of cohomology of line
bundles on $X$ is an almost numerical condition, in the sense that
in most cases we can decide whether or not a particular cohomology group
vanishes  by simply looking at the numerical class to which
the line bundle belongs.
As a matter of fact, in those numerical classes in which we
cannot decide,  there exist  line bundles for
which certain cohomology group vanishes and  line bundles for which
it does not.  We will study in more detail this situation in the
case
$e = -1$, because we will need   for the sequel
 to know exactly for which line
bundles the cohomology vanishes. Concretely,  this knowledge will
allow us to  use Proposition 2.1 and Proposition 5.4  in the proofs
of Theorem 4.2 and Theorem 5.7
respectively. It will be used as well in [GP].
Also we
will show the existence of a smooth elliptic curve numerically
equivalent to
$2C_0-f$.
\proclaim {Proposition 3.2} Let $X$ be a ruled surface with
invariant $e = -1$. Then

$3.2.1.$ There exist only three effective line bundles in the
numerical class of $2C_0 -f$. They are $\Cal O (2C_0 - (\frak e +
\eta _i)f)$, where the $\eta _i$s are the nontrivial degree $0$
divisors corresponding to the three nonzero torsion points
in Pic$^0(C)$. The unique element in $|2C_0 - (\frak e + \eta
_i)f)|$ is  a smooth elliptic curve $E_i$

$3.2.2.$ For each $n > 1$, there are only four effective line
bundles numerically equivalent to $n(2C_0 -f)$. They are $\Cal O
(2nC_0 - n(\frak e + \eta
_i)f)$ and $\Cal O (2nC_0 - n\frak e f)$. The only smooth
(elliptic) curves (and indeed the only irreducible curves) in these
numerical classes are general members in $|4C_0 - 2\frak e f|$.

The number of linearly independent global sections of these line
bundles is
summarized in the following table:
\vskip .2 cm
$$\vbox {\offinterlineskip \hrule
\halign{\vrule
#&\hfil\quad#\quad\hfil&\vrule#&\hfil\quad#\quad\hfil&\vrule#&
\hfil\quad#\quad\hfil&\vrule#&\hfil\quad#\quad\hfil&\vrule#&
\hfil\quad#\quad\hfil&\vrule#&
\hfil\ #\ \hfil&\vrule#&
\hfil\quad#\quad\hfil&\vrule#
\cr
height .1cm &&height .1cm &&height .1cm &&height .1cm&&height
.1cm&&height .1cm&&height
.1cm&&height .1cm
  \cr
  &$n \geq 0$&&$0$&&$1$&&$2$&&$3$&&&&$n$& \cr
height .1cm &&height .1cm &&height .1cm &&height .1cm&&height .1cm
&&height .1cm &&height .1cm&&height .1cm \cr
\noalign{\hrule}
height .1cm &&height .1cm &&height .1cm &&height .1cm&&height .1cm
  &&height .1cm&&height .1cm&&height .1cm\cr
  &$\text h ^0(\Cal O (2nC_0 - n \frak e
f))$&&$1$&&$0$&&$2$&&$1$&&&&$3\lfloor\frac {n} {2}\rfloor-n+1$& \cr
height .1cm &&height .1cm &&height .1cm &&height .1cm&&height .1cm
  &&height .1cm &&height .1cm&&height .1cm\cr
\noalign{\hrule}
height .1cm &&height .1cm &&height .1cm &&height .1cm&&height .1cm
  &&height .1cm &&height .1cm&&height .1cm\cr
  &$\text h ^0(\Cal O (2nC_0 - n (\frak e + \eta _i)
f))$&&$0$&&$1$&&$1$&&$2$&&&&$n-\lfloor\frac {n} {2}\rfloor$& \cr
height .1cm &&height .1cm &&height .1cm &&height .1cm&&height .1cm
  &&height .1cm&&height .1cm&&height .1cm\cr }
\hrule }$$

\endproclaim
\vskip .2 cm

{\it Proof.}
For any $p \in C$ we consider the following exact sequence:
$$\displaylines{0 \to \text H ^0(\Cal O (2C_0 - pf)) \to \text H
^0(\Cal O (2C_0)) @> \varphi _p >> \text H ^0(\Cal O _{pf} (2C_0))
\cr
\to \text H ^1(\Cal O (2C_0-pf)) \to 0 \ .\cr}$$
Pushing forward the morphism
$$\text H^0(\Cal O(2C_0)) \otimes \Cal O \to \Cal O(2C_0)$$
to $C$ we obtain
$$ \text H ^0(S^2(\Cal E)) \otimes \Cal O _C @> \varphi >> S^2(\Cal E) \to
Q \to 0 \ .$$ Note that the restriction of $\varphi$ to the fiber
of $\text H ^0(S^2(\Cal E))\otimes \Cal O _C$ over $p$   is
precisely $\varphi_p$. Thus the points
$p$ for which H$^0(\Cal O (2C_0 -pf))
\neq 0$ are exactly the ones where the rank of $\varphi$ drops.
Note that  h$^0(\Cal O (2C_0))=3$ (push down the bundle to $C$, use
the same semistability considerations as in the sketch of the
proof of Proposition 3.1 to obtain the vanishing of H$^1$ and then,
use Riemann-Roch.) The rank of
$S^2(\Cal E)$ is also 3, so if the rank of
$\varphi$ never dropped, $\varphi$ would be an isomorphism, which
is not true, because, since $e=-1$, the degree of $S^2(\Cal E)$ is
$3$. Therefore   there exists
$p
\in C$ such that H$^0(\Cal O (2C_0 -pf)) \neq 0$. We fix such a
point $p$ and some
 effective divisor $E$ inside $|2C_0-pf|$. Since $2C_0-f$
cannot be written as sum of two nonzero numerical classes both
containing  effective divisors, $E$ is irreducible
and reduced. By adjunction, $p_a(E) =1$ and since $E$ dominates
$C$, $E$ is indeed a smooth elliptic curve.

We prove now by
induction on
$n $ the following statement: for each \text{$n\geq 0$}, there are
finitely many effective line bundles numerically equivalent to
\linebreak
\text{$2nC_0 -nf$}. The result is obviously true for
$n=0$. Take now $n>0$. We fix  a divisor $\frak d '$ of degree $n-1$.
What we want to prove is that the number of points $z \in C$ such
that H$^0(\Cal O (2nC_0 -\frak d f)) \neq 0$ is finite , where
$\frak d=
\frak d ' + z$. We may assume that H$^0(\Cal O
((2n-2)C_0 -(\frak d ' -p +z) f)=0$, since, by induction
hypothesis, there are only finitely many points $z$ for which this
does not happen. We tensor the sequence
$$0 \to \Cal O(-2C_0+pf) \to \Cal O \to \Cal O _E \to
0 \eqno(3.2.3)$$ by $\Cal O (2nC_0 -\frak d f)$ and take global
sections. Since H$^0(\Cal O
((2n-2)C_0 -(\frak d  -p) f)$ is $0$ and
the degree of the push forward of $\Cal O((2n-2)C_0 -(\frak d
-p)f)$ to $C$  is 0, it follows that H$^1(\Cal O
((2n-2)C_0 -(\frak d ' -p +z) f)=0$. Hence we obtain that
$\text H^0(\Cal O(2nC_0 -\frak d f) = \text H ^0(\Cal
O_E(2nC_0-\frak d f))$. The degree of $\Cal O_E(2nC_0-\frak d f)$
is zero, therefore
$z$ is such that H$^0(\Cal O(2nC_0
-\frak d f) \neq 0$ iff $\Cal O_E(2nC_0-\frak d f)
\simeq \Cal O _E$, i.e., iff $\Cal O_E(2nC_0-\frak d '
f)
\simeq \Cal O _E (zf)$. There are only finitely many such points
$z$, since otherwise, we will have that all the fibers of the
$2:1$ morphism from $E$ onto $C$ induced by the degree 2
divisor which is obtained as the restriction of $2nC_0-\frak d '
f$ to $E$ are members of the same
$g^1_2$.

The last statement implies that the length
of
$Q$ is finite, and equal to  \linebreak \text{$\text{deg}\ S^2(\Cal E) = 3$.}
We claim
that
$Q$ is in fact supported in three distinct points $p_1, p_2,$ \linebreak
$p_3$. If not, there would exist a global section of $S^2(\Cal E)$
vanishing at some
\text{point $q$} to order greater or equal than $2$. In particular,
$\Cal O (2C_0 -2qf)$ would be effective, which contradicts
Proposition 3.1. Our aim now is to identify $p_1, p_2, p_3$. Let
$E_i$ be the unique element of
$|2C_0 - p_if|$. We saw before that $E_i$ is a smooth elliptic
curve. Pushing down to $C$ the exact sequence (3.2.3) we obtain
$$0 \to \Cal O _C \to \pi _* \Cal O _{E_i} \to  R ^1\pi_* \Cal
O (-2C_0 +p_if) \to 0 \ .$$
Since $\pi |_{E_i}$ is unramified and $E_i$ is connected, it
follows that $\pi_*\Cal O_{E_i} = \Cal O_C \oplus \Cal L$ for some line
bundle $\Cal L$ such that
$\Cal L ^{\otimes 2} = \Cal O_C$, but $\Cal L \neq \Cal O_C$. Using
relative duality
and projection formula one obtains that
$2(p_i-\frak e) \sim 0$ but $p_i-\frak e \not\sim
0$. This proves the first part of the proposition.

For the second part, remember that we have already proven the
existence of only finitely many effective line bundles.
Therefore, for any $p \in C$, we have the exact sequence
$$ 0 \to \text H ^0(S^{2n}\Cal E \otimes \Cal O_C ((-n+1)p)) \otimes
\Cal O_C \to  S^{2n}\Cal E \otimes \Cal O_C ((-n+1)p) \to Q' \to 0 \
.$$
The length of $Q'$ is equal to the degree of $S^{2n}\Cal E \otimes \Cal
O_C ((-n+1)p)$, which is $4m +1$ if $n = 2m$ and $4m +3$ if $n=
2m +1$, but it is also equal to the sum of the dimensions of the
linear spaces of global sections of line bundles in the numerical
class of $2nC_0 -nf$. Then, the rest of the statement in 3.2.2 and
the numbers in the table  follow from comparing the length of
$Q'$ with the sum of the dimensions of the linear spaces generated
by sections corresponding to reducible divisors numerically
equivalent to
$2nC_0 -nf$. \boxit{}

\vskip .2 cm
(3.2.4) We will fix once and for all a smooth elliptic curve $E$ in
the numerical class of $2C_0-f$.
\vskip .2 cm
(3.2.5) For a different proof of the existence of a smooth elliptic
curve in the numerical class of $2C_0-f$ see [Ho2], corollary
{2.2.}
\vskip .2 cm
In the case $e \geq 0$ we are interested in finding  sections of
$\pi$ whose self-intersection is near to that of $C_0$ (they will
play in the sequel a role similar to that of
$E$):

\proclaim {Proposition 3.3} Let $X$ be an elliptic ruled surface
with invariant $e \geq 0$. The general member of
$|C_0 -\frak e f|$ is a smooth elliptic curve and
those are the only smooth curves in the numerical
class of $C_0 +ef$.
\endproclaim

{\it Proof.}
If det$( \Cal E) \neq \Cal O$
the dimension
of
$|C_0-\frak e f|$ is $\text h^0(\Cal O \oplus \Cal O(-\frak e))
-1 = e$. Since the dimension of $|-\frak e' f|$ is $e-1$ for any
nontrivial divisor
$-\frak e'$ of degree $e$ on $C$
it is clear that not all the elements in $|C_0-\frak e f|$ are
unions of $C_0$ and $e$ fibers. On the other hand this is the only
way in which an element of $|C_0-\frak e f|$ can be
reducible (this is because for any divisor $\frak d$ of
degree $d < e$, the dimension of $|C_0 +
\frak df|$ is $d-1$, which implies that any element of $|C_0 +
\frak df|$ is the union of $C_0$ and $d$ fibers). Thus the
general member of $|C_0-\frak e f|$ is irreducible. Therefore, it
maps surjectively onto $C$ and hence it is a smooth elliptic curve. If
$\frak e' \equiv \frak e$ but $\frak e' \neq \frak e$, then
$$\text{dim}|C_0-\frak e' f|  = \text h^0(\Cal O (\frak e -\frak e')
\oplus
\Cal  O(-\frak e')) -1 =
e -1$$
which means that all members of $|C_0-\frak e' f|$ are
reducible.

If det$( \Cal E) = \Cal O$, $ \Cal E$ is an extension of $\Cal
O$ by $\Cal O$. Thus the member or members of $|C_0-\frak e f| =
|C_0|$ are smooth elliptic curves and $|C_0-\frak e' f| =
\emptyset$ for any divisor $\frak e'$ on $C$ of degree $0$
different from $\frak e$.
\boxit{}

\vskip .2 cm
(3.3.1) We will fix once and for all a smooth elliptic curve $E'$ in
the numerical class of $C_0 + ef$.

\proclaim {Proposition 3.4 ([H], V.2.20.b and V.2.21.b.)}
Let $L$ be line bundle on $X$ in the numerical class of $aC_0 +
bf$.

If $e =-1$, $L$ is ample iff $a >0$ and $a >-\frac 1 2 b$.

If $e \geq 0$, $L$ is ample iff $a >0$ and $b-ae >0$.
\endproclaim

\vskip .4 cm
We state now a proposition describing numerical conditions which
imply base-point-freeness. The proof follows basically the
one given in [Ho1] and [Ho2]. In characteristic $0$
the
proposition can also be proven using Reider's theorem ([R]).

\proclaim {Proposition 3.5}
Let $L$ be a line bundle on $X$ in the numerical class of $aC_0 +
bf$.

If $e =-1$ , $a \geq 0$,
$a+b \geq 2$ and 	$a+2b \geq 2$, then $L$ is base-point-free.

If $e \geq 0$, $a \geq 0$ and
$b-ae \geq 2$, then $L$ is base-point-free.
\endproclaim

{\it Proof.}
First we consider the case $e = -1$. In the first place we prove
the proposition when $a =0$ and $b \geq 2$, when $a=b=1$, and
when $a=2$ and $b=0$. The first case follows easily from the fact
that line bundles on elliptic curves whose degrees are greater or
equal than $2$ are base-point-free. For the  other two cases we use
the fact that  there are only one or two
minimal sections through a given point of $X$ (c.f. [Ho2]).
On the other
hand, for a given point $p \in C$,  there are infinitely many
effective reducible divisors in
$|C_0 + pf|$, namely, those consisting of the union of a divisor
linearly equivalent to $C_0 + \tau f$ and a divisor linearly
equivalent to $(p - \tau)f$, where $\tau$ is a degree $0$ divisor
on $C$. Hence, the intersection of all those reducible divisors is
empty. Analogously, for a given divisor $\nu$ of degree $0$
  there are infinitely many  effective reducible divisors in
$|2C_0 +  \nu f|$, namely, those consisting of the union of a
divisor linearly equivalent to $C_0 + (\nu +\tau) f$ and a divisor
linearly equivalent to $C_0 + (\nu - \tau)f$, and the same argument goes
through.

Now we use lemma 2.4. The base-point-free line bundle $B$ will be
numerically equivalent to $bf \ (b \geq 2)$, $C_0 +f$ or  $2C_0$ and
$Y$ will be $E$ (defined in (3.2.4)) or
$C_0$ (note that since $\text{deg}(L\otimes \Cal O_Y(Y))
\geq 2$, it follows that
$L \otimes \Cal O_Y(Y)$ is base-point-free). Iterating this process
we obtain the result. The only place where we have to be careful in
the application of lemma 2.4 iteratively is in making sure that the
base-point-free line bundles we keep obtaining are nonspecial. This
problem is taken care of  by proposition 3.1.

The case $e \geq 0$ is easier. The line bundle $L$ is
base-point-free if it is in the numerical class of $bf$,  when $b
\geq 2$. Then we get the result for any other bundle satisfying the
conditions in the proposition by using the lemma 2.4. The curve
$Y$ in lemma 2.4 will be  $E'$ (defined in (3.3.1)).  Again
proposition 3.1 assures us  that the line bundles we obtain are
nonspecial. \boxit{}

\vskip .3 cm

\proclaim{Remark 3.5.1} The numerical condition of Proposition 3.5
characterizes those equivalence classes consisting entirely of
base-point-free line bundles.
\endproclaim

{\it Proof.}   A line bundle that satisfies the above numerical
conditions is base-point-free by virtue of Proposition 3.5. To prove the other
implication, consider a base-point-free line bundle
$L$ in the numerical class of $aC_0+bf$, which does not satisfy the
above conditions. If $e(X)=-1$, the restriction of $L$ to the elliptic curve
$E$ is a base-point-free line bundle. Hence, since its degree is equal to
$a+2b < 2$, it must be the trivial line bundle, which implies that $a+2b=0$.
Then if follows from
Proposition 3.2 that the general member of the numerical class
is not  base-point-free (in fact it is not even effective!).
If $e(X) \geq 0$,  for the same reason as above, the restriction of $L$ to
$C_0$
is trivial. This is only possible if $L = \Cal O(n(C_0-\frak e f))$.\boxit{}

(3.5.2) The proof of the previous remark suggests that there  exist nontrivial
base-point-free line bundles with self-intersection $0$. That is indeed the
case.
For example, if $e = -1$, the divisors $2nC_0 +n \frak e
f$  for any even number $n$ greater than
$0$ are base-point-free;  if $X = C
\times \bold P ^1$, the divisors $nC_0$   and if
$e \geq 1$, the divisors $n(C_0 - \frak e f)$ are base-point-free.
Hence base-point-freeness cannot be characterized numerically.
However, if we assume $L$ to be ample, then the numerical
conditions in Proposition 3.5 do give a characterization of
base-point-freeness:

\proclaim {Remark 3.5.3}
Let $L$ be a line bundle on $X$ in the numerical class of $aC_0 +
bf$.

If $e =-1$ , the line bundle $L$ is ample and base-point-free iff $a \geq 1$,
$a+b \geq 2$ and 	$a+2b \geq 2$.

If $e \geq 0$, the line bundle $L$ is ample and base-point-free iff $a \geq 1$
and
$b-ae \geq 2$.
\endproclaim

{\it Proof.}
If $L$ satisfies the numerical conditions in the statement of the remark,
then by propositions 3.4 and 3.5 it is ample and base-point-free.
Now assume that $L$ is ample and base-point-free. If $e
=-1$, from proposition 3.4, it follows that $a \geq 1$. On
the other hand, since $L$ is base-point-free, its restriction to any curve in
$X$ is also base-point-free. Consider the curves $C_0$ and a smooth curve
$E$ in the numerical class of $2C_0-f$. The restriction of $L$ to each of
them has degree $a+b$ and $a+2b$ respectively. The fact that the restriction
of $L$ to $C_0$ is base-point-free implies that either $a+b \geq 2$ or the
restriction of $L$ to $C_0$ is trivial. The latter is impossible since $L$ is
ample. Analogously the fact that $L$ is ample and that the restriction of $L$
to $E$ is base-point-free implies that $a+2b \geq 2$. If $e \geq 0$, by
proposition 3.4, $a \geq 1$. Since $L$ is as well base-point-free,
by restricting $L$ to $C_0$ we obtain that $b-ae = \text{deg} \ (L
\otimes \Cal O_{C_0}) \geq 2$.
\boxit{}

\heading 4. Normal presentation on elliptic ruled surfaces
\endheading

We recall that in this section $X$ denotes a ruled surface over
an elliptic curve and we continue to use the notation introduced
at the beginning of section 3. We have just seen which line bundles
on
$X$ are ample and which are base-point-free. The next question to
ask would be: ``which line bundles are very ample and which are
normally generated?". This problem was solved by Y. Homma in [Ho1]
and [Ho2], who proved that a line bundle $L$ on $X$ is normally
generated  iff it is very ample. She also characterizes
those line bundles (see Figures 1 for the case $e=-1$).
Homma  proves as well that in the case of a normally generated
line bundle $L$,  the ideal corresponding to the embedding induced
by
$L$ is generated by quadratic  and cubic forms. Thus the next
question is to identify those line bundles which are normally
presented.
\vskip .2 cm

(4.1) Throughout the remaining part of this section we assume that
\text{$\text{char}(\bold k)
\neq 2$}.
\vskip .2 cm
We will use Proposition 2.1 and the results from
Section 3 to characterize the line bundles on $X$ which are
normally presented.

\proclaim{Theorem 4.2}
The condition of normal presentation depends only on numerical
equivalence. More precisely, let
$L$ be a line bundle on
$X$ numerically equivalent to
$aC_0+bf$. If $e = -1$, $L$ is normally presented iff
$a \geq 1$, $a+b \geq 4$ and $a+2b \geq 4$. If $e \geq 0$, $L$ is
normally presented iff $a \geq 1$ and $b-ae \geq 4$.
\endproclaim

{\it Proof.}
First we prove that if a line bundle $L$ satisfies the numerical
conditions in the statement, it is normally presented. To this end
we will use Proposition 2.1. The idea is to write $L$ as tensor
product of two line bundles $B_1$ and $B_2$ satisfying the
numerical conditions in Proposition 3.5, and such that  H$^2(B_1
\otimes B_2^*)=\text H ^2(B_2
\otimes B_1^*)=0$.
Let us exhibit  the line bundles $B_1$ and $B_2$ in the
different cases:

If $e=-1$, $L$ can be written as tensor product of $B_1$ and
$B_2$, where the couple $B_1$ and $B_2$ satisfies one of the
following numerical conditions:
\vskip .1 cm
$(4.2.1) \qquad B_1 \equiv C_0+nf$ and $B_2 \equiv 2f$ or $C_0 +f$,
for some $n
\geq 1$; (in this case, $1 \leq a \leq 2$ and $a+b \geq 4$).
\vskip .05 cm
$(4.2.2) \qquad B_1 \equiv 2C_0$ and $B_2 \equiv 2C_0 +lf$
or
$ C_0+nf
$, for some $l \geq 0$, and some $n
\geq 1$; (in this case $3 \leq a \leq 4$ and $a+b \geq 4$).
\vskip .05 cm
$(4.2.3) \qquad B_1 \equiv 2C_0 + m(2C_0 -f)$ and $B_2 \equiv 2C_0
+lf$ or $ C_0+nf
$, for some $m \geq 1$,
some $l \geq 1$ and some $n
\geq 1$; (in this case, $a \geq 5$ and $a +2b
>4$).
\vskip .05 cm
$(4.2.4) \qquad B_1 \equiv 2C_0 + m(2C_0 -f)$ and $B_2 \equiv
2C_0$, for some $m \geq 1$; (in this case, $a \geq 5$ and $a +2b
=4$).
\vskip .1 cm
If $B_1$ and $B_2$ satisfy (4.2.1), (4.2.2) or (4.2.3),
Proposition 3.1 implies that \linebreak
\text{$\text  H^2(B_1
\otimes B_2^*)=\text H ^2(B_2
\otimes B_1^*)=0$.} If $B_1$ and $B_2$ satisfy $(4.2.4)$,
H$^2(B_1
\otimes B_2^*)$ vanishes but $\text H ^2(B_2
\otimes B_1^*)$ might not be zero. However from Proposition
3.2 it follows that, given $L \equiv 4C_0 + m(2C_0 -f)$, $m \geq
1$, we can choose $B_1 \equiv 2C_0 + m(2C_0 -f)$ and $B_2 \equiv
2C_0$ such that $L = B_1 \otimes B_2$ and $\text H ^2(B_2
\otimes B_1^*)=0$, hence we are done in this case.

If $e \geq 0$, $L$ can be written as tensor product of $B_1$
and
$B_2$, where,
\vskip .1 cm
if $a$ is even: $B_1 \equiv (a/2) C_0 + \lfloor b/2 \rfloor f$ and
$B_2
\equiv (a/2) C_0 +
\lceil b/2 \rceil f$ \quad \ \ and
\vskip .05 cm
 if $a$ is odd: $B_1 \equiv \lfloor a/2 \rfloor C_0 +
\lfloor (b-e)/2 \rfloor f$ and
$B_2
\equiv \lceil a/2
\rceil C_0 +
\lceil (b+e)/2 \rceil f$.
\vskip .4 cm
Proposition 3.4 implies that $L$ is
ample, Proposition 3.5 implies  that $B_1$ and $B_2$ are
base-point-free and Propositions 3.1 and 3.2 imply that H$^1(B_1)
= \text H ^1(B_2)=
\text H ^2(B_1
\otimes B_2^*)=\text H ^2(B_2
\otimes B_1^*)=0$, so from Proposition 2.1 it follows that $L$ is
normally presented.
\vskip .2 cm

Now we will suppose that $L$ is  normally
presented but does not satisfy the
numerical conditions in the statement and we will derive
a contradiction. We can assume that $a \geq 1$. Otherwise $L$ would not be
ample.   If
$e = -1$, we can also assume that
$a +b
= 3$ and $a+2b = 3$. If not the restriction of $L$ to either the
minimal section $C_0$ or  the curve $E$ defined in (3.2.4) would not
be very ample. Analogously, if
$e
\geq 0$,  we can assume that  $b -ae = 3$. Otherwise, the
restriction of $L$ to
$C_0$ would not be very ample.

To obtain the
contradiction we follow the same strategy: we will see that the
assumption of $L$ being normally presented forces its restriction
to both $C_0$ and $E$  to be also
normally presented, and we will derive from that the
contradiction. If
$e
\geq 0$ let
$P$ denote the line bundle $
\Cal O (C_0)$. If $e=-1$ and $L \equiv C_0 + 2f$ or
$2C_0 +f$, let
$P$ denote
$\Cal O (C_0)$.  If $e = -1$ and $L
\equiv 3C_0 + m(2C_0 -f)$,  $m
\geq 0$, let $P$ denote $ \Cal O (E)$. Let $\frak p$ denote a
smooth elliptic curve in
$|P|$. We claim that
$$ \text H^2(\bigwedge
^2M_L
\otimes L
\otimes P ^*)=0 \ . \eqno (4.2.5)$$  To
prove (4.2.5) we will prove instead the fact that
\text {H$^2(M_L^{\otimes 2}
\otimes L
\otimes P ^*)=0$ .} Consider the following exact sequence, which
arises from exact sequence (1.1),
$$\displaylines{\hfill \text H ^1(M_L \otimes
L ^{\otimes 2}
\otimes P ^*) \to \text H ^2(M_L^{\otimes 2} \otimes L
\otimes P ^*) \hfill\cr
\hfill \to \text H ^0(L) \otimes \text H ^2 (M_L
 \otimes L
\otimes P ^*)\ . \hfill\llap{(4.2.6)}\cr}$$
Since $\text H ^1(L) = 0$, the vanishing of   $\text H ^1(M_L
\otimes L ^{\otimes 2}
\otimes P ^*)$ is equivalent to the surjectivity of
the  map
$$\text H^0(L) \otimes \text H^0(L^{\otimes 2} \otimes P^*) \to \text
H^0(L^{\otimes 3}\otimes P^*)\ .$$
The line bundle
$L$ is base-point-free by Proposition 3.5. Therefore, by Theorem
1.3,
it is enough to check that $\text H^1(L
\otimes P ^*)= \text H ^2(P ^*)=0$. The vanishings follow
from our choice of
$P$ and from Proposition 3.1, except the vanishing of H$^2( P
^*)$ when
$P
\simeq
\Cal O (E)$, which follows from Proposition 3.2 and duality.
Using (1.1) we obtain that H$^2(M_L
\otimes L
\otimes P ^*)$ will vanish if $\text H^1(L^{\otimes 2} \otimes P^*)$ and
$\text H^2(L \otimes P^*)$ vanish. These two vanishings follow from
Proposition 3.1. Therefore by (4.2.6), H$^2(M_L^{\otimes 2}
\otimes L
\otimes P ^*)=0$ and H$^2(\bigwedge
^2M_L \otimes L
\otimes P ^*)=0$.

Now since $L$ is assumed to be normally
presented and H$^1(L) = 0$ we have, by Lemma 1.2, that
H$^1(\bigwedge ^2M_L \otimes L )=0$. Thus from (4.2.5) it follows
that
$$\text H
^1(\bigwedge
^2(M_L \otimes \Cal O_\frak p) \otimes L  ) =
\text H^1(\bigwedge ^2M_L \otimes L  \otimes \Cal
O_\frak p)= 0 \ . \eqno (4.2.7)$$
Consider the
following commutative diagram, (which holds for any
base-point-free line bundle $L$  and for any nontrivial line bundle
$P$ such that \linebreak
$\text H^1(L \otimes P ^*)=0$):
$$\matrix
&&0&&0&& \\
&& \downarrow && \downarrow && \\
0 &\to& \text H ^0(L \otimes P ^*) \otimes \Cal O _\frak p &
\to & \text H ^0(L \otimes P ^*) \otimes \Cal O _\frak p
&\to& 0 \\
&& \downarrow && \downarrow && \downarrow \\
0& \to & M_L \otimes \Cal O _\frak p & \to  & \text H ^0(L
) \otimes \Cal O_\frak p
&\to& L \otimes \Cal O_\frak p & \to & 0 \\
&& \downarrow && \downarrow && \downarrow \\
0& \to& M_{L \otimes \Cal O _\frak p}&\to&\text H ^0(L \otimes
\Cal O _\frak p) \otimes \Cal O _\frak p& \to &L \otimes
\Cal O _\frak p &\to &0 \\
&& \downarrow && \downarrow && \downarrow \\
&&0&&0&&0 \\
\endmatrix$$
{}From the left hand side vertical sequence we obtain the surjection $$
\bigwedge
^2(M_L \otimes \Cal O_\frak p) \otimes L
\to
\bigwedge
^2M_{L \otimes \Cal O_\frak p} \otimes L  \ .$$
Since $\frak p$ is a curve, it follows that
$$\text H
^1(\bigwedge
^2M_{L \otimes \Cal O_\frak p} \otimes L  ) = 0 \ . \eqno (4.2.8)$$ The line
bundle $L \otimes \Cal O_\frak p$ on $\frak p$  is normally
generated because
$\frak p$ is an elliptic curve and $\text{deg}\ (L \otimes \Cal
O_\frak p) = 3$. Thus Lemma 1.2 and (4.2.8) imply that   $L \otimes
\Cal O_\frak p$ is normally presented, which is impossible since the complete
linear series of $L
\otimes
\Cal O_\frak p$ embeds $\frak p$ as a plane cubic! \boxit{}

\vskip .2 cm

(4.2.9) It follows in particular from Proposition 4.2 that the
normally generated line bundles on $X$, or more precisely, their
numerical classes form a convex set in $\text{Num}(X)$, as shown
in Figure 1, in which we describe $\text{Num}(X)$ when $e=-1$ (we do not
draw the similar picture for $e \geq 0$).

\vskip .2 cm
Now we reformulate Theorem 4.2 and state some corollaries, which will help us
to
put our result in perspective:

\proclaim{Theorem 4.3} Let
$X$ be an elliptic ruled surface. A line
bundle
$L$ on
$X$ is normally presented iff it is ample and can be written as the
tensor product of two line bundles
$B_1$ and
$B_2$ such that every line bundle numerically equivalent to any of
them is base-point-free.
\endproclaim

{\it Proof.} If $L$ is normally presented, it is obviously ample
and satisfies the numerical conditions of Proposition 4.2. In the
proof of that proposition we showed that a line bundle satisfying
the mentioned numerical conditions can be written as the tensor
product of two base-point-free line bundles $B_1$ and $B_2$
satisfying the conditions of Remark 3.5.1. Hence these $B_1$ and
$B_2$ are such that all the line bundles in their numerical classes
are base-point-free.

On the other hand, assume that $L$ is ample and isomorphic to $B_1
\otimes B_2$ and that any line bundle numerically equivalent to
either
$B_1$ or $B_2$ is base-point-free. Let
$B_i$ be in the numerical class of
$a_iC_0 + b_if$ and $L$ in the numerical
class of $aC_0+bf$. If $e=-1$, by Remark 3.5.1,
$a_i +
b_i \geq 2$ and
$a_i+2b_i
\geq 2$. Thus we obtain that
$
a+b \geq 4$ and  $a+2b \geq 4$.
Since $L$ is ample, $a\geq 1$. Hence, from Proposition 4.2 it
follows that $L$ is normally presented. If $e \geq 0$ one argues in a similar
fashion.
\boxit{}

\proclaim{Corollary 4.4} Let $X$ be as above. Let $B_i$ be an ample
and base point free line bundle on $X$ for all $1 \leq i \leq q$.
If $q \geq 2$, then
$B_1 \otimes \dots \otimes B_q$ is normally presented and if $q <2$, in
general $B_1 \otimes \dots \otimes B_q$ is not normally presented.
\endproclaim

{\it Proof.}  It follows from Remarks 3.5.1 and 3.5.3
that a line bundle numerically equivalent to any of the $B_i$ is
base-point-free. From Remark 3.5.3 and Proposition 3.4 it follows
that $L$ is ample. Thus, by Theorem 4.3, $L$ is normally
presented.
\qed\par
\epsfverbosetrue
\epsfbox{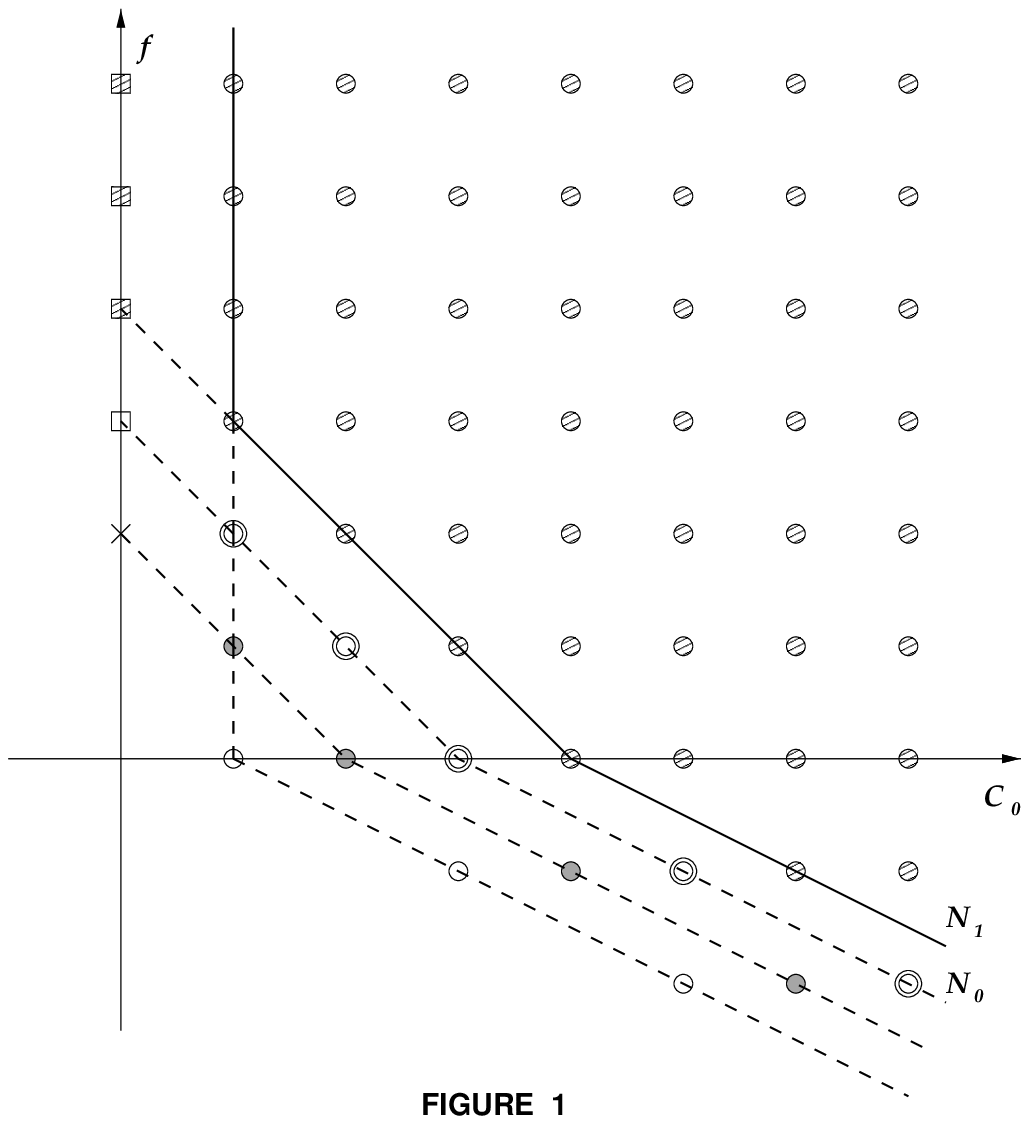}

In the above figure cross means that all the members in the numerical
class are base point free, dashed (or lined) square means that the
corresponding coordinate ring is presented
by quadratic forms, dashed (or lined) disc means normally
presented, annulus means normally
generated, blank disc means ample, gray or hashed disc means ample
and base point free.

\proclaim{Corollary 4.5} Let $X$ be as above. Let $A_i$ be an ample
line bundle on $X$ for all $1 \leq i \leq q$. If $q \geq 4$, then
$A_1 \otimes \dots \otimes A_q$ is normally presented and  if $q <4$, in
general $A_1 \otimes \dots \otimes A_q$ is not normally presented.
\endproclaim

{\it Proof.} From Proposition 3.4 and Remark 3.5.3 it follows that
the tensor product of two ample line bundles is ample and
base-point-free. Hence the corollary follows from Corollary 4.4.
\boxit{}
\vskip .2 cm

\proclaim{Corollary 4.6 } Let
$X$ be an elliptic ruled surface. Let $A_i$ be an
ample line bundle on
$X$ for all $1
\leq i \leq q$.

If
$e=-1$ and $q \geq 5$, then $\omega_X \otimes A_1 \otimes \dots \otimes A_q$
is normally presented. If $e=-1$ and $q < 5$, in general
$\omega_X \otimes A_1 \otimes \dots \otimes A_q$
is not normally presented.

If
$e= 0$ and $q \geq 4$, then $\omega_X \otimes A_1 \otimes \dots
\otimes A_q$ is normally presented. If $e=0$ and $q < 4$, in
general
$\omega_X \otimes A_1 \otimes \dots \otimes A_q$
is not normally presented.

If $e \geq 1$ and $q \geq 3$, then $\omega_X \otimes A_1 \otimes
\dots \otimes A_q$ is normally presented. If $e \geq 1$ and $q <
3$, in general $\omega_X \otimes A_1 \otimes \dots \otimes A_q$ is
not normally presented.
\endproclaim

{\it Proof.}
Let $A_i$ be in the numerical class of $a_iC_0 +b_if$
and $\omega _X \otimes A_1 \otimes \dots \otimes A_q$
in the numerical class of $aC_0+bf$. If $e=-1$, $A_i$ is ample
iff $a_i \geq 1$ and $a_i +2b_i \geq 1$ (c.f. Proposition 3.4). In
particular we also have that if $A_i$ is ample, then
$a_i+b_i \geq 1$. Since $\omega _X $ is numerically
equivalent to $-2C_0+f$ it follows that
$$\displaylines {a\geq q-2 \geq 3 > 1 \cr
a+b \geq q-1 \geq 4 \quad \text{and} \cr
a+2b \geq q \geq 4 \ .}$$
Hence by Theorem 4.2, $\omega _X \otimes A_1 \otimes \dots
\otimes A_q$ is normally presented.

If $e = 0$, $A_i$ is ample
iff $a_i \geq 1$ and $b_i \geq 1$ (c.f. Proposition 3.4).
Since
$\omega _X$ is numerically
equivalent to $-2C_0$ it follows that $a \geq q-2 >1$
and $b-ae \geq q \geq 4$. Hence by Theorem 4.2, $\omega _X \otimes
A_1
\otimes \dots
\otimes A_q$ is normally presented.

If $e \geq 0$, $A_i$ is ample
iff $a \geq 1$ and $b_i - a_ie \geq 1$ (c.f. Proposition 3.4).
Since
$\omega _X $ is numerically
equivalent to $-2C_0-ef$ it follows that $a \geq q-2 \geq 1$
and $b-ae \geq q+e \geq 4$. Then by Theorem 4.2, $\omega _X \otimes
A_1 \otimes \dots
\otimes A_q$ is normally presented.

The line bundles $\omega_X \otimes \Cal O(4C_0)$, if $e(X)=-1$;
$\omega_X
\otimes
\Cal O(3C_0 +3f)$, if $e(X)=0$; and $\omega_X \otimes \Cal
O(2C_0+2(e+1)f)$, if $e(X) \geq 1$, are not normally presented (c.f.
Theorem 4.2). Thus our bound is sharp.
\boxit{}

\vskip .4 cm

We want to compare our results to the results known for curves. Fujita and
St. Donat (c.f. [F] and [S-D] ) proved that, on a curve, any line bundle of
degree bigger or
equal than $2g+2$ is normally presented. The results in this section as well as
Theorem 2.1 are analogous in different ways to the result by Fujita and St.
Donat.
The approach taken up to now to
generalize this result has been to look at adjoint linear series.
In this line
it was conjecture by Mukai that on any surface $X$, $\omega_X \otimes A_1
\otimes \dots \otimes A_q$
should be normally presented for
all $q \geq 5$ and $A_i$ ample line bundle. Corollary 4.6 shows that this
conjecture holds if $X$ is an elliptic ruled surface and that it is sharp if
the invariant $e(X)=-1$.
\pageno=17
One disadvantage of this generalization is that it only gives information about
 a
small class of line bundles. The possible ways of generalization indicated by
Theorem 4.4:

(4.7) Let $X$ be a surface.
If $L$ is the product of two ample line bundles $B_1$ and $B_2$, such that
every line bundle $B$ numerically equivalent to either $B_1$ or $B_2$ is
base-point-free, then $L$ is normally presented;

or by Proposition 2.1:

(4.8) Let $X$ be a surface.
If $L$ is ample and the product of two base-point-free and nonspecial
 line bundles, then $L$ is normally presented;

or maybe by some combination of the two, take in account a larger class of line
bundles in general. In subsequent articles we prove that both (4.7) and (4.8)
hold for K3 surfaces.

 We remark that Theorem 4.3, which is stronger than
Corollary 4.4, can also be seen as an analogue of Fujita and St. Donat's
theorem,
since
the latter
can be rephrased as follows:

Let $L$ be a line bundle on a curve. Every line bundle numerically equivalent
to
$L$ is normally presented iff $L$ is ample and the tensor product of two line
bundles $B_1$ and $B_2$ such that every line bundle $B$ numerically equivalent
to $B_1$ or $B_2$ is base-point-free.

However, the veracity of Theorem 4.3 seems to depend on the particular
properties
of elliptic ruled surfaces and the corresponding statement is false on K3
surfaces.
\vskip .2 cm
We generalize Corollaries 4.4 and 4.5 to higher syzygies in a forthcoming
article
(c.f. [GP]), by proving that the product of $p+1$ or more ample and
base-point-free
line bundles satisfies the property $N_p$.

\heading 5. Koszul algebras
\endheading

In the previous section we determined which line bundles on an
elliptic ruled surface are normally presented. A question to ask
is whether the coordinate ring of the embedding induced by those
line bundles is a Koszul ring, since it is well known
that a variety with a Koszul
homogeneous coordinate ring is projectively normal and defined by
quadrics.
The answer to this question is affirmative not only in
the  case of elliptic ruled surfaces, but in all other cases with
which we
have dealt throughout this work, since we are able to prove that
the corresponding  coordinate ring to a line bundle
satisfying the conditions of Proposition 2.1 is Koszul.

We introduce now some  notation and some basic
definitions:
given a line bundle
$L$ on a variety
$X$, we recall that $R(L) = \bigoplus _{n=0} ^\infty \text H
^0(X,L ^{\otimes n})$.

\proclaim {Definition 5.1} Let $R= \bold k \oplus R_1 \oplus R_2
\oplus
\dots$ be a graded ring and $\bold k$ a field. $R$ is a Koszul ring
iff Tor$_i^R(\bold k,\bold k)$
has pure degree $i$ for all $i$.
\endproclaim

Now we will give a cohomological interpretation, due to
Lazarsfeld, of the Koszul
property for a coordinate ring of type $R(L)$. Let $L$ be a
globally  generated line
bundle on a variety
$X$. We will denote $M^{(0),L} := L$ and $M^{(1),L} :=M_L \otimes
L = M_{M^{(0),L}} \otimes L$. If $M^{(1),L}$ is globally generated,
we denote
$M^{(2),L}  :=M_{M^{(1),L}} \otimes L$.  We repeat the process
and define inductively $M^{(h),L}  :=M_{M^{(h-1),L}} \otimes L$,
if $M^{(h-1),L}$ is
globally generated. Now we are ready to state
the following

\proclaim {Lemma 5.2 ([P], Lemma 1)}
Let $X$ be a projective variety over an algebraic closed field
$\bold k$. Assume that $L$ is a base-point-free line bundle on $X$
such that the vector bundles $M^{(h),L}$ are globally generated for
every $h \geq 0$.
If \text{H$^1(M^{(h),L} \otimes L^{\otimes s}) =0$} for
every $h \geq 0$ and every $s \geq 0$ then $R(L)$
is a Koszul $\bold k$-algebra. Moreover, if
H$^1(L^{\otimes s})=0$ for every $s \geq 1$ the converse is also
true.
\endproclaim

Now we will prove a general result analogous to Proposition 2.1 but
before that, we state the following well known

\proclaim {Observation 5.3}
Let $\Cal F$ be a locally free sheaf over a scheme $X$ and $A$ an
ample line bundle.
If the multiplication map H$^0( \Cal F \otimes A ^{\otimes n})
\otimes
\text H ^0(A) \to \text H ^0( \Cal F \otimes A ^{\otimes n+1})$
surjects for all $n \geq 0$, then $\Cal F$ is globally generated.
\endproclaim

\proclaim {Theorem 5.4}
Let $X$ be a surface with $p_g = 0$, let $B_1$ and $B_2$ be two
base-point-free line bundles and let $L = B_1 \otimes B_2$ be
ample. If
$$\text H^1(B_1)=\text H ^1(B_2)=\text H ^2(B_1
\otimes B_2 ^*) = \text H ^2(B_2 \otimes B_1 ^*)=0 \ ,$$
then the following properties are satisfied for all $h \geq 0$:
\vskip .3 cm
$5.4.1)$ \ $M^{(h),L}$ is globally generated
\vskip .05 cm
$5.4.2)$ \ H$^1(M^{(h),L} \otimes B_1 ^{\otimes b_1} \otimes B_2
^{\otimes b_2})=0$ for all $ b_1, b_2 \geq 0$
\vskip .05 cm

$5.4.3)$ \ H$^1(M^{(h),L} \otimes B_j
^*)=0$ where $j =1,2$
\vskip .05 cm
$5.4.4)$ \ H$^1(M^{(h),L}\otimes B_i \otimes B_j
^*)=0$ where $i =1,2$ and $j=2,1$
\vskip .05 cm
$5.4.5)$ \ H$^1(M^{(h),L}\otimes B_i ^{\otimes 2} \otimes B_j
^*)=0$ where $i =1,2$ and $j=2,1$
\vskip .3 cm
In particular H$^1(M^{(h),L} \otimes L
^{\otimes s})=0$ for all $h,s \geq 0$, and $R(L)$ is a Koszul
$k$-algebra.
\endproclaim

{\it Proof.}
We prove the lemma by induction on $h$.

If $h=0$, property 5.4.1 means that $L$ is globally generated,
which is true by hypothesis. Properties 5.4.2 to 5.4.5 mean that
H$^1(B_1 ^{\otimes (b_1+1)} \otimes B_2
^{\otimes (b_2+1)}) = \text H^1(B_i ^{\otimes \beta _i}) =0$ where
$b_1, b_2 \geq 0$,  $\beta _i = 1,2,3$ and $i=1,2$. These
vanishings occur by hypothesis and Observation 2.2.

Now consider $h >0$ and assume that the result is true for all $0
\leq h' \leq h-1$. Let $L'$ denote $B_1 ^{\otimes b_1} \otimes B_2
^{\otimes b_2}$, $B_j
^*$, $B_i \otimes B_j
^*$ or $ B_i ^{\otimes 2} \otimes B_j
^*$ accordingly. If we tensor
$$0 \to M_{M^{(h-1),L}} \to \text H^0(M^{(h-1),L})\otimes \Cal O
\to M^{(h-1),L} \to 0\ , \eqno (\text{c.f.}\  (1.1))$$
by $L \otimes L'$ and take global sections, we obtain
$$\displaylines{\text
H^0(M^{(h-1),L})
\otimes
\text H ^0(L \otimes L') @>\alpha>> \text H ^0(M^{(h-1),L}\otimes
L \otimes L') \to \cr
\hfill\text H^1(M^{(h),L} \otimes L') \to \text
H^0(M^{(h-1),L})
\otimes
\text H ^1(L \otimes L')\ .\hfill\llap{(5.4.6)}}$$ Since,
by Observation 2.2,
$\text H^1(L
\otimes L') =0$,
properties 5.4.2
to  5.4.5 are equivalent to the surjectivity of the multiplication
map $\alpha$ in the different cases.
First we prove property 5.4.5. Consider the following commutative diagram:
$$\matrix
\text H^0(M^{(h-1),L})
\otimes
\text H ^0(B_i) \otimes
\text H ^0(B_i) \otimes
\text H ^0(B_i) \ \to&  \text H^0(M^{(h-1),L})
\otimes \text H ^0(B_i^{\otimes 3})\cr
@VV\varphi_1 V\cr
\text H^0(M^{(h-1),L} \otimes B_i)
\otimes
\text H ^0(B_i) \otimes
\text H ^0(B_i)&@VV\alpha V\cr
@VV\varphi_2 V \cr
\text H^0(M^{(h-1),L} \otimes B_i^{\otimes 2})
\otimes
\text H ^0(B_i) \qquad \ @>\varphi_3>>&\text H^0(M^{(h-1),L}
 \otimes B_i ^{\otimes 3}) \ .\cr
\endmatrix$$

To show the surjectivity of $\alpha$ it suffices then to show
the surjectivity of $\varphi_1$, $\varphi_2$ and $\varphi_3$.
To prove that these three map are surjective we
use Theorem 1.3. For example, to see that $\varphi_1$ is surjective is enough
by Theorem 1.3 to
show that $\text H ^1(M^{(h-1),L} \otimes B_i
^*) =0$ and $\text H ^2(M^{(h-1),L} \otimes B_i ^{-2}) =0$. We argue
analogously
for the other two maps and thus conclude that in order to show the surjectivity
of $\alpha$ it is enough to check that
$$\displaylines{\rlap{(5.4.7)}\hfill \text H ^1(M^{(h-1),L} \otimes B_i
^*) =
\text H ^1(M^{(h-1),L}) = 0 \hfill \cr
\hfill \text H ^1(M^{(h-1),L} \otimes B_i) = 0 \quad \text{and}\hfill \cr
\rlap{(5.4.8)} \hfill \text H ^2(M^{(h-1),L} \otimes B_i ^{-2}) = \text
H ^2(M^{(h-1),L} \otimes B_i ^*) = 0 \hfill \cr
\hfill \text H ^2(M^{(h-1),L} ) =
0 \ .\hfill \cr}$$
\vskip .2 cm

The vanishings in (5.4.7) follow from the assumption that properties
5.4.2 and
5.4.3 hold for
$h - 1$. Proving  (5.4.8) is not hard. For instance, to obtain
\vskip .3 cm
$(5.4.8.1) \qquad
\text H ^2(M^{(h-1),L} \otimes B_i ^{-2}) =0$
\vskip .2 cm
we consider the  following
sequence that we obtain from (1.1):
$$\displaylines{\text H^1(M^{(h-2),L} \otimes B_j \otimes B_i ^{*}
) \to \cr \text
H^2(M^{(h-1),L}
\otimes B_i ^{-2})
 \to
\text H^2( B_j
\otimes B_i ^{*})\ .}$$
Hence it is clear that in order to show (5.4.8.1) it is enough to check that
\linebreak
$
\text H ^1(M^{(h-2),L} \otimes B_j
\otimes B_i ^*) =0$ and that $\text H ^2(B_j  \otimes B_i
^*) =0$. Arguing in a similar way for the remaining vanishings in (5.4.8),
we conclude that in order to prove (5.4.8) it is enough to check
\vskip .05 cm
$$\displaylines{\rlap{(5.4.9)} \hfill \text H ^1(M^{(h-2),L} \otimes B_j
\otimes B_i ^*) = \text H ^1(M^{(h-2),L} \otimes B_j)=0 \hfill \cr
\hfill\text H
^1(M^{(h-2),L} \otimes B_1 \otimes B_2) = 0 \quad \text{and} \hfill \cr
\rlap{(5.4.10)} \hfill \text H ^2(B_j  \otimes B_i
^*) = \text H
^2(B_j) = \text H ^2(B_1 \otimes B_2) = 0 \ . \hfill \cr}$$
\vskip .2 cm

The vanishings in (5.4.9) follow from the assumption that properties 5.4.2
and 5.4.4 hold for
$h-2$. Statement (5.4.10) follows by hypothesis and
Observation 2.3.

The proof of properties  5.4.3 and 5.4.4 is analogous. In fact, notice that we
have implicitly proven both when we showed the surjectivity of
$\varphi_1$ and $\varphi_2$.

Now we prove property 5.4.2. The argument is similar to the one we
have use to prove 5.4.5 and we will only sketch it here in little
detail.  To show the surjectivity of the map $\alpha$
$$  \text H^0(M^{(h-1),L})
\otimes
\text H ^0(B_1^{\otimes b_1+1} \otimes B_2^{\otimes b_2+1}) \to
\text H^0(M^{(h-1),L} \otimes B_1 ^{\otimes b_1 +1} \otimes B_2
^{\otimes b_2 +2})$$
(c.f. (5.4.6)), one can write a similar diagram to the one in the proof of
5.4.5.  Then it is  enough to prove the surjectivity of the following map,
which is a composition of multiplication maps (we assume $b_2 \geq b_1$):
$$ \matrix \text H^0(M^{(h-1),L})
\otimes
[\text H ^0(B_1) \otimes \text H ^0(B_2)]^{\otimes b_1 + 1}
\otimes
[\text H ^0(B_2)]^{\otimes b_2 - b_1} \cr   @VV\varphi V \cr
 \text H^0(M^{(h-1),L} \otimes B_1 ^{\otimes b_1 +1} \otimes B_2
^{\otimes b_2 +2})  \cr \endmatrix$$

We show the surjectivity of the composite map $\varphi$ by showing the
surjectivity of each of its components. The first component is
$$ \matrix \text H^0(M^{(h-1),L})
\otimes \text H ^0(B_1) \otimes \text H ^0(B_2) \otimes
[\text H ^0(B_1) \otimes \text H ^0(B_2)]^{\otimes b_1}
\otimes
[\text H ^0(B_2)]^{\otimes b_2 - b_1} \cr   @VV\varphi_1V \cr
 \text H^0(M^{(h-1),L} \otimes B_1) \otimes \text H ^0(B_2) \otimes
[\text H ^0(B_1) \otimes \text H ^0(B_2)]^{\otimes b_1}
\otimes
[\text H ^0(B_2)]^{\otimes b_2 - b_1} \cr \endmatrix$$

Hence by Theorem 1.3
it is enough to check the vanishings of the cohomology groups
\text{$\text H^1(M^{(h-1),L}\otimes B_1^*)$}
and $\text H^2(M^{(h-1),L}\otimes B_1^{-2})$. For the surjectivity of
the
second component
$$ \matrix \text H^0(M^{(h-1),L} \otimes B_1) \otimes \text H ^0(B_2)
\otimes
[\text H ^0(B_1) \otimes \text H ^0(B_2)]^{\otimes b_1}
\otimes
[\text H ^0(B_2)]^{\otimes b_2 - b_1} \cr
   @VV\varphi_2V \cr
 \text H^0(M^{(h-1),L} \otimes B_1  \otimes B_2) \otimes
[\text H ^0(B_1) \otimes \text H ^0(B_2)]^{\otimes b_1}
\otimes
[\text H ^0(B_2)]^{\otimes b_2 - b_1} \ ,\cr \endmatrix$$
again by Theorem 1.3 it
is enough to check the  vanishings of the groups
\linebreak
$\text H^1(M^{(h-1),L}\otimes B_1
\otimes B_2^*)$
and $\text H^2(M^{(h-1),L}\otimes B_1 \otimes  B_2^{-2})$.
We use the same argument for the remaining components of $\varphi$ and conclude
that in order to prove the surjectivity of $\varphi$, it suffices to check
\vskip .3 cm
$(5.4.11) \qquad \text H^1(M^{(h-1),L} \otimes B_1 ^{\otimes \beta
_1 }
\otimes B_2 ^{\otimes \beta _2 }) =0$, for all $\beta_1$ and
$\beta_2$ satisfying one of the following conditions:
\vskip .1 cm
(5.4.11.1) \qquad $-1
\leq
\beta _1
\leq b_1 -1$ and $\beta _2 = \beta _1 +1$
\vskip .05 cm
(5.4.11.2) \qquad $1 \leq \beta _1 \leq
b_1 +1$ and $\beta _2 = \beta _1 -2$
\vskip .05 cm
(5.4.11.3)\qquad $\beta _1 = b_1
+1$ and
$b_1 \leq \beta _2 \leq b_2 -1$
\vskip .2 cm
and
\vskip .2 cm
$(5.4.12) \qquad \text H^2(M^{(h-1),L} \otimes B_1 ^{\otimes \gamma
_1 }
\otimes B_2 ^{\otimes \gamma _2 }) =0$,  for all $\gamma
_1$ and $\gamma
_2$ satisfying one of the following conditions:
\vskip .1 cm
(5.4.12.1) \qquad $-2 \leq \gamma
_1
\leq
b_1 -2$ and $\gamma _2 = \gamma _1 +2$
\vskip .05 cm
(5.4.12.2) \qquad $1 \leq \gamma _1
\leq b_1
+1$ and $\gamma _2 = \gamma _1 -3$
\vskip .05 cm
(5.4.12.3) \qquad $\gamma _1 = b_1
+1$ and $b_1 -1\leq \gamma _2 \leq b_2 -2$ \.
\vskip .2 cm

The vanishings in $(5.4.11)$, except the vanishings of $\text
H^1(M^{(h-1),L}
\otimes B_1 ^{* } )$ and $\text H^1(M^{(h-1),L}
\otimes B_1 \otimes B_2^*)$,
follow from the assumption  that property 5.4.2 holds  for $h-1$.
The
vanishing of
$\text H^1(M^{(h-1),L}
\otimes B_1 ^{ * } ) $ follows from the assumption  that
property 5.4.3 holds for $h-1$. The
vanishing of $\text H^1(M^{(h-1),L} \otimes B_1 \otimes B_2^*)$
follows from the assumption  that
property 5.4.4 holds for $h-1$.
To prove the vanishings in $(5.4.12)$ we consider the  following
sequence that we obtain from (1.1):
$$\displaylines{\text H^1(M^{(h-2),L} \otimes B_1 ^{\otimes (\gamma
_1 + 1)}
\otimes B_2 ^{\otimes (\gamma _2 +1) }) \to \cr \text
H^2(M^{(h-1),L}
\otimes B_1 ^{\otimes \gamma _1 }
\otimes B_2 ^{\otimes \gamma _2 }) \to
\text H^2( B_1 ^{\otimes
(\gamma _1 + 1)}
\otimes B_2 ^{\otimes (\gamma _2 +1) })\ .}$$
Hence it is enough to show
that these cohomology groups vanish:

\vskip .3 cm
$(5.4.13) \qquad \text H^1(M^{(h-2),L} \otimes B_1 ^{\otimes (\gamma
_1 + 1)}
\otimes B_2 ^{\otimes (\gamma _2 +1) }) =0$ \qquad and
\vskip .05 cm
$(5.4.14) \qquad \text H^2( B_1 ^{\otimes (\gamma _1
+ 1)}
\otimes B_2 ^{\otimes (\gamma _2 +1) }) =0$,
\vskip .05 cm
for all
$\gamma _1$ and $\gamma _2$ satisfying one of the conditions from
(5.4.7.1) to (5.4.7.3).
\vskip .2 cm
Statement $(5.4.13)$, except for the vanishings of $\text
H^1(M^{(h-2),L} \otimes B_2
\otimes B_1 ^{ -1 } )$ and $\text
H^1(M^{(h-2),L} \otimes B_1^{\otimes 2}
\otimes B_2 ^{ -1 } )$, follow from the assumption
 that property 5.4.2  holds  for $h-2$. The vanishing of $\text
H^1(M^{(h-2),L} \otimes B_2
\otimes B_1 ^{ -1 }
)$ follows from the assumption that property 5.4.4 holds for
$h-2$. The vanishing of $\text
H^1(M^{(h-2),L} \otimes B_1^{\otimes 2}
\otimes B_2 ^{ -1 } )$ follows from the assumption that property 5.4.5 holds
for
$h-2$. All the vanishings in
$(5.4.14)$ follow by hypothesis and Observation 2.3.

\vskip .1 cm

Finally we prove property 5.4.1. By
Observation 5.3, it is enough to show that the map
$$ \text H ^0(M^{(h),L} \otimes L ^{\otimes n})
\otimes
\text H ^0(L ) \to \text H ^0(M^{(h),L} \otimes L ^{\otimes n+1}) \eqno
(5.4.15)
$$ surjects for all $n \geq 0$.  For that it suffices to prove
the
surjectivity of the map
$$ \text H ^0(M^{(h),L} \otimes L ^{\otimes n}) \otimes \text H
^0(B_1) \otimes \text H
^0(B_2)\to \text H ^0(M^{(h),L} \otimes L ^{\otimes n+1}) $$
for all $n \geq 0$. Using Theorem 1.3, it is enough to check

$$\displaylines{\rlap{(5.4.16)}\hfill\text H ^1(M^{(h),L} \otimes B_1
^{\otimes n-1}
\otimes B_2 ^{\otimes n}) =0\hfill\cr
\hfill\text H ^1(M^{(h),L} \otimes B_1
^{\otimes n+1}
\otimes B_2 ^{\otimes n-1}) = 0 \hfill \llap{\text{and}\quad} \cr
\rlap{(5.4.17)}\hfill\text H ^2(M^{(h),L} \otimes B_1 ^{\otimes
n-2}
\otimes B_2 ^{\otimes n}) =0\hfill\cr
\hfill\text H ^2(M^{(h),L} \otimes B_1
^{\otimes n+1}
\otimes B_2 ^{\otimes n-2}) = 0 \hfill\cr}$$
\vskip .1 cm
The vanishings in (5.4.16) follow from the fact, which we have just proved,
that properties 5.4.2 to 5.4.4 hold for $h$. To prove (5.4.17),
again by (1.1), it is enough to show that
$$\displaylines{\rlap{(5.4.18)} \hfill
\text H ^1(M^{(h-1),L} \otimes B_1
^{\otimes n-1}
\otimes B_2 ^{\otimes n+1}) = 0\hfill\cr
\hfill
\text H ^1(M^{(h-1),L} \otimes
B_1 ^{\otimes n+2}
\otimes B_2 ^{\otimes n-1}) = 0\hfill\llap{\text{and}\qquad}\cr
\rlap{(5.4.19)} \hfill
\quad \text H ^2(B_1 ^{\otimes n-1}
\otimes B_2 ^{\otimes n+1}) = \text H ^2(B_1
^{\otimes n+2}
\otimes B_2 ^{\otimes n-1}) = 0 \ . \hfill \cr}$$

\vskip .1 cm
The vanishings in (5.4.18) follow from the assumption that
properties  5.4.2, 5.4.4 and 5.4.5 hold for
$h-1$ and (5.4.19) follows by hypothesis and
Observation 2.3.
\vskip .1 cm
In particular,  it follows from  property 5.4.2
that  H$^1(M^{(h),L} \otimes L
^{\otimes s})=0$ for all $h,s \geq 0$. Thus, as consequence of
Lemma 5.2, the coordinate ring $R(L)$ is a Koszul
$\bold k$-algebra. \boxit{}
\vskip .35 cm
(5.5) Note that if $h=1$ and $n=1$, the multiplication map (5.4.15)
is actually the same  as (2.7.1). Moreover, the fact that
\text{$\text H^1(M_L \otimes L) =0$} is a special case of 5.4.2. Hence
on our way
to
prove Theorem 5.4, we have reproved Proposition 2.1 and therefore Theorem
 5.4 may be seen as a generalization of the cited proposition.
 \vskip .2 cm
Even though the above theorem is stated for surfaces with $p_g=0$, the same
proof (or indeed a simpler one) works for curves. Thus we obtain the following

\proclaim {Theorem 5.5}
Let $C$ be a curve, let $B_1$ and let $B_2$ be two nontrivial
base-point-free line bundles
on $C$. If $\text H^1(B_1)=\text H^1(B_2)=0$,
then $R(L)$ is a Koszul.
\endproclaim
{\it Proof.}
The only properties of surfaces
with $p_g=0$ that we use in the proof of Theorem 5.4 are the fact that $\text
H^2
(\Cal O_X)=0$, and Observations 2.2 and 2.3. Observation 2.2 is obviously still
true if $X$ is a curve. Observation 2.3 and the fact that $\text H^2
(\Cal O_X)=\text H^2(B_1^*\otimes B_2)=\text H^2(B_2^*\otimes B_1)=0$
and are trivially true for curves, hence the theorem follows from the proof of
Theorem 5.4. \boxit{}
\vskip .2 cm
Theorem 5.5 yields as a corollary the following result by David Butler (see
also [Po]):
\proclaim{Corollary 5.6 ([B], Theorem 3)} Let $C$ be a curve and let $L$ be a
line bundle on $C$. If $\text{deg}(L) \geq 2g+2$, then $R(L)$ is Koszul.
\endproclaim
{\it Proof.}
If $\text{deg}(L) \geq 2g+2$, then $L$ can be written as tensor product of two
general line bundles of degree $g+1$. Such line bundles are base-point-free and
nonspecial. \boxit{}
\vskip .2 cm
Theorem 5.4 yields these three results:

\proclaim {Corollary 5.7}
Let $X$ be an Enriques surface over an algebraic closed field of
characteristic $0$ and let
$B$ be an ample line bundle on $X$ without base points. Then
$R(B^{\otimes 2})$ is Koszul.
\endproclaim

{\it Proof.}
The proof is analogous to the proof of Corollary 2.8. \boxit{}

\proclaim{Theorem 5.8}
Let $X$ be an elliptic ruled surface. Let $L$ be a line bundle on $X$
numerically equivalent to
$aC_0+bf$. If $e = -1$
 and $a \geq 1$, $a+b \geq 4$ and $a+2b \geq 4$, then
$R(L)$ is Koszul. If $e \geq 0$
and $a \geq 1$ and $b-ae\geq
4$, then $R(L)$ is Koszul.
\endproclaim

{\it Proof.}
The proof is analogous to the proof of the first part of Theorem
4.2. \boxit{}
\vskip .2 cm
(5.9)
It is well known that for an ample line bundle $L$, the fact of
$R(L)$ being Koszul implies formally the
property of being normally presented (c.f. [BF], 1.16).
Therefore Theorem 5.5 gives a different proof of Fujita and St. Donat's
theorem,
Corollary 5.7
provides another proof of Corollary 2.8 and Theorem 5.8
provides another proof of the first part of Theorem 4.2. These
proofs are less elementary, but in the case of Theorem 4.2, we have
the advantage of working also in characteristic 2.
\vskip .3 cm
If we
assume that $\text{char}(\bold k)\neq 2$, it follows from Proposition 4.2
that the property of $R(L)$ being a Koszul algebra is characterized by the
numerical conditions in the statement of Theorem 5.8. We can restate this
as we did in the case of Theorem 4.2:

\proclaim{Theorem 5.9} Let
$X$ be as above and let $L$ be a line bundle on $X$. Assume that
$\text{char}(\bold k)\neq 2$. Then $R(L)$ is a Koszul algebra
iff it is ample and can be written as the tensor product of two line
bundles
$B_1$ and
$B_2$ such that every line bundle numerically equivalent to any of
them is base-point-free.
\endproclaim

(5.9.1)  Having in account that, on elliptic ruled surfaces,
normal presentation only depends on numerical equivalence (c.f. Theorem 4.2),
Theorem 5.9 can be
considered as analogous to Theorem 5.5.
Indeed Theorem 5.5  can be rephrased  as follows:

If a line bundle $L$ on $C$ is normally presented and so is every line
bundle
numerically equivalent to $L$, then $R(L)$ is Koszul.

\heading  references \endheading

\roster

\item"[BF]" J. Backelin \& R. Fr\"oberg, {\it Veronese subrings,
and rings with linear resolutions}, Rev. Roum. Math. Pures Appl.
{\bf 30} (1985), 85-97.

\item
"[B]" D. Butler, {\it Normal generation of vector bundles over
a curve}, J. Differential
Geometry {\bf 39} (1994) 1-34.

\item"[C]" G. Castelnuovo, {\it Sui multipli di uni serie di
gruppi di punti appatente ad una
curva algebraica}, Rend. Circ. Mat.
Palermo {\bf 7} (1892) 99-119.

\item"[F]" T. Fujita, {\it Defining equations for certain types of polarized
varieties}, Complex Analysis and Algebraic Geometry, Cambridge Univ. Press,
1977,
165-173.

\item"[GP]" F.J. Gallego \& B.P. Purnaprajna, {\it Higher syzygies of elliptic
ruled surfaces}, Preprint.

\item"[G]" M. Green, {\it Koszul cohomology and the geometry of
projective varieties},
J. Differential Geometry {\bf
19} (1984) 125-171.

\item"[GL]" M. Green \& R. Lazarsfeld, {\it Some results on the
syzygies of finite sets and
algebraic curves}, Compositio
Math. {\bf 67} (1989) 301-314.

\item"[H]" R. Hartshorne {\it Algebraic Geometry}, Springer,
Berlin, 1977.

\item"[Ho1]" Y. Homma, {\it Projective normality and the
defining equations of ample
invertible sheaves on elliptic ruled
surfaces with $e \geq 0$}, Natural Science
Report, Ochanomizu
Univ. {\bf 31} (1980) 61-73.

\item"[Ho2]"
\hbox{\leaders \hrule  \hskip .6 cm}\hskip .05 cm , {\it Projective
normality and the defining equations of an elliptic ruled
surface with negative
invariant}, Natural Science Report, Ochanomizu
Univ. {\bf 33} (1982) 17-26.

\item"[Mi]" Y. Miyaoka, {\it The Chern class and Kodaira
dimension of a minimal variety},
Algebraic Geometry --Sendai 1985,
Advanced Studies in Pure Math.,
Vol. 10,
North-Holland,
Amsterdam, 449-476.

\item "[Mu]" D. Mumford,
{\it Varieties defined by quadratic equations}, Corso CIME in
Questions on Algebraic Varieties,
Rome, 1970, 30-100.

\item "[P]" G. Pareschi, {\it Koszul algebras associated to
adjunction bundles}, J. of Algebra
{\bf 157} (1993) 161-169.

\item "[Po]" A. Polishchuk, {\it On the Koszul property of the
Homogeneous coordinate ring of a curve}, J of Algebra
{\bf 178} (1995) 122-135

\item"[R]" I. Reider, {\it Vector bundles of rank $2$  and
linear systems on an algebraic surface},
Ann. of Math. (2) {\bf
127} (1988) 309-316.

\item "[S-D]" B. St.-Donat, {\it Sur les equations d\'efinissant
 une courbe alg\'ebrique, C.R. Acad. Sci. Paris, Ser. A {\bf 274} (1972),
324-327.

\endroster

\enddocument